\documentclass[11pt,a4paper,reqno]{amsart}

\usepackage[T1]{fontenc}
\usepackage[utf8]{inputenc}
\usepackage{lmodern}
\usepackage[margin=25mm]{geometry}
\usepackage{amsmath,amssymb,amsthm,mathtools}
\usepackage{microtype}
\usepackage{enumitem}
\usepackage{cite}
\usepackage[hidelinks]{hyperref}

\hypersetup{
  pdftitle={Smoothing of the Coupled KdV System on the Torus},
  pdfauthor={Tugrul Burak Gurel and Murat Keçkin},
  pdfsubject={Analysis of partial differential equations},
  pdfkeywords={Coupled KdV system, Boussinesq equations, nonlinear smoothing}
}

\newcommand{\R}{\mathbb{R}}
\newcommand{\Z}{\mathbb{Z}}
\newcommand{\T}{\mathbb{T}}
\newcommand{\wtr}{\widetilde{R}}

\theoremstyle{plain}
\newtheorem{theorem}{Theorem}[section]
\newtheorem{lemma}[theorem]{Lemma}
\newtheorem{proposition}[theorem]{Proposition}

\theoremstyle{definition}

\theoremstyle{remark}
\newtheorem{remark}[theorem]{Remark}
\numberwithin{equation}{section}

\allowdisplaybreaks[4]
\setlist[itemize]{leftmargin=*,topsep=5pt,itemsep=4pt,parsep=0pt}
\date{}

\begin{document}

\title[Smoothing of the Coupled KdV System on $\mathbb{T}$]{Smoothing of the Coupled KdV System on $\mathbb{T}$}
\author[Tu\u{g}rul Burak G\"urel \and Murat Ke\c{c}kin]{Tu\u{g}rul Burak G\"urel \and Murat Ke\c{c}kin}

\address{Tu\u{g}rul Burak G\"urel \\ Department of Mathematics, Bo\u gazi\c ci University, 
Bebek 34342, Istanbul, Türkiye}
\email{bgurel@boun.edu.tr}

\address{Murat Ke\c{c}kin \\ Department of Mathematics, Bo\u gazi\c ci University, 
Bebek 34342, Istanbul, Türkiye}
\email{murat.keckin@std.bogazici.edu.tr}

\subjclass[2020]{}


\begin{abstract}
We study the nonlinear smoothing properties of the Coupled KdV system posed on the torus. This system arises as a particular case of the four-parameter family of Boussinesq systems studied by Bona, Chen and Saut, corresponding to the choice $\eta^{2}=\frac{2}{3}$, $b=d=0$ and $a=c=\frac{1}{6}$. Although the similarities with the periodic KdV equation and the favorable resonance structure of the Majda-Biello system for the coupling parameter $\alpha =-1$ might suggest that smoothing can be established by a relatively straightforward adaptation of existing techniques, the nonlinear interactions in the Coupled KdV system introduce additional difficulties. In particular, certain mixed interactions lead to degeneracies in the associated resonance functions, resulting in a loss of frequency decay that restricts the attainable smoothing gain. By combining normal-form reductions, duality arguments and Bourgain-space estimates, we establish a smoothing gain of any order
$
a< \min \{1,\frac{2}{3}(2s+1)\}
$
for $s>-\frac{1}{2}$. Our result demonstrates that the interplay between nonlinear interactions and resonance structures plays a decisive role in determining the smoothing properties of coupled dispersive systems.
\end{abstract}

\maketitle

\section{Introduction}
\label{sec:intro}
In this paper, we investigate the nonlinear smoothing properties of the Coupled KdV system
\begin{equation}
\label{eq:coupled-kdv}
\begin{cases}
    u_{t}+v_{x}+(uv)_{x}+\dfrac{1}{6}v_{xxx}=0,\\[2mm]
    v_{t}+u_{x}+vv_{x}+\dfrac{1}{6}u_{xxx}=0,\\[2mm]
    (u,v)|_{t=0}=(u_{0},v_{0}),
\end{cases}
\end{equation}
posed on the periodic domain $\mathbb{T}\times \mathbb{R}$, where $\mathbb{T}=\mathbb{R}/(2\pi \mathbb{Z})$.

Before turning to the Boussinesq equations, we briefly review the principal developments in the well-posedness theory for the Korteweg-de Vries (KdV) equation
\begin{equation}
\begin{cases}
        u_{t}+u_{xxx}+uu_{x}=0,\\
        u(x,0)=u_{0}(x)\in H^{s}(\mathbb{A}),
    \end{cases}
    \qquad \mathbb{A}=\mathbb{R} \text{ or } \mathbb{T}.
\end{equation}
A major advance in the study of the KdV equation was made by Bourgain, who introduced Fourier restriction spaces and thereby extended local well-posedness to the $L^{2}$ setting on both $\mathbb{R}$ and $\mathbb{T}$ \cite{Bourgain-partII}. Subsequently, Kenig, Ponce and Vega established local well-posedness for initial data in $H^{s}(\mathbb{R})$ when $s>-\frac{3}{4}$ and in $H^{s}(\mathbb{T})$ when $s>-\frac{1}{2}$ \cite{KenigPonceVega1996}. The endpoint case $s=-\frac{3}{4}$ on the real line was later resolved by Christ, Colliander and Tao \cite{ChristCollianderTao2003}. In the periodic setting, Colliander, Keel, Staffilani, Takaoka and Tao proved both local and global well-posedness for data in $H^{-\frac{1}{2}}(\mathbb{T})$ \cite{CollianderEtAl2003}. At the endpoint $s=-\frac{3}{4}$ on $\mathbb{R}$, global well-posedness was independently established by Guo \cite{Guo2009} and Kishimoto \cite{Kishimoto2009}. Exploiting the integrable structure of the KdV equation, Kappeler and Topalov obtained local and global well-posedness in $H^{-1}(\mathbb{T})$ \cite{KappelerTopalov2006}. In contrast to these well-posedness results, Molinet demonstrated that the equation is ill-posed in $H^{s}(\mathbb{A})$ for $\mathbb{A}=\mathbb{R}$ or $\mathbb{A}=\mathbb{T}$ whenever $s<-1$ \cite{Molinet2011,Molinet2012}. Most recently, Killip and Vi\c{s}an established global well-posedness in $H^{-1}(\mathbb{R})$ \cite{KillipVisan2019}. These developments have substantially advanced the understanding of the KdV equation at low regularity, bringing its well-posedness theory to a particularly comprehensive stage.

The system \eqref{eq:coupled-kdv} belongs to the class of Boussinesq-type systems, which arise in the modeling of the bidirectional propagation of long, small-amplitude gravity waves on the surface of water. These systems describe the interplay between nonlinear wave interactions and linear dispersion, both of which play an essential role in the evolution of long waves. The classical Boussinesq equations, originally introduced to describe the propagation of surface waves in a canal, provide a fundamental example of this class; see \cite{Boussinesq1871,Boussinesq1872}. In the original physical formulation, $u+h_{0}$ represents the total water depth, whereas $v$ denotes the horizontal fluid velocity, with $h_{0}$ being the undisturbed water depth. A general four-parameter family of Boussinesq-type systems is given by
\begin{equation}
\label{eq:abcd-system}
\begin{cases}
      u_{t}+v_{x}+(uv)_{x}+av_{xxx}-bu_{xxt}=0,\\
      v_{t}+u_{x}+vv_{x}+cu_{xxx}-dv_{xxt}=0,\\
      (u,v)|_{t=0}=(u_{0},v_{0}),
    \end{cases}
\end{equation}
where the coefficients satisfy
\begin{equation*}
a+b=\frac{1}{2}\left(\eta^{2}-\frac{1}{3}\right), \qquad c+d=\frac{1}{2}\left(1-\eta^{2}\right), \qquad a+b+c+d=\frac{1}{3}.
\end{equation*}
Different choices of the parameters $(a,b,c,d)$ lead to various models within this family, including the Kaup system, the Bona-Smith system and the classical Boussinesq system. In particular, the Coupled KdV system \eqref{eq:coupled-kdv} considered in this work corresponds to the choice $\eta^{2}=\frac{2}{3}$, $a=c=\frac{1}{6}$ and $b=d=0$. Bona, Chen and Saut \cite{Bona-Saut-Chen-Media-I,Bona-Saut-Chen-M-II} developed a systematic theory for this family of models, including their derivation from the two-dimensional Euler equations governing free-surface flows and the analysis of their well-posedness.

To analyze \eqref{eq:coupled-kdv}, we first diagonalize its linear part. Let 
\begin{equation*}
U=\begin{pmatrix}
    u \\v
\end{pmatrix}, \qquad W=P^{-1}U, \qquad P=\begin{pmatrix}
    1 &-1 \\ 1&1
\end{pmatrix}.
\end{equation*}
After relabeling the components of $W$ and their initial data as $(u,v)$ and $(u_{0},v_{0})$, respectively, the system becomes
\begin{equation}
\label{eq:oldCoupled-to-oldcoupled}
\begin{cases}
        u_{t}+u_{x}+\frac{1}{6}u_{xxx}+\frac{3}{4}(u^{2})_{x}-\frac{1}{4}(v^{2})_{x}+\frac{1}{2}(uv)_{x}=0,\\
        v_{t}-v_{x}-\frac{1}{6}v_{xxx}-\frac{1}{4}(u^{2})_{x}+\frac{3}{4}(v^{2})_{x}+\frac{1}{2}(uv)_{x}=0,\\
        (u,v)|_{t=0}=(u_{0},v_{0}), \qquad (x,t)\in \mathbb{T}\times \mathbb{R}.
    \end{cases}
\end{equation}
We then introduce the changes of variables
\begin{equation*}
\widetilde{u}(x,t)=u(x+t,t), \qquad \widetilde{v}(x,t)=v(x-t,t).
\end{equation*}
After dropping the tildes, we obtain the system
\begin{equation}
\label{eq:coupled-to-coupled}
\begin{cases}
        u_{t}+\frac{1}{6}u_{xxx}+\frac{3}{4}(u^{2})_{x}-\frac{1}{4}(v^{2})_{x}+\frac{1}{2}(uv)_{x}=0,\\
        v_{t}-\frac{1}{6}v_{xxx}-\frac{1}{4}(u^{2})_{x}+\frac{3}{4}(v^{2})_{x}+\frac{1}{2}(uv)_{x}=0,\\
        (u,v)|_{t=0}=(u_{0},v_{0}), \qquad (x,t)\in \mathbb{T}\times \mathbb{R}.
    \end{cases}
\end{equation}
Throughout the remainder of the paper, $u$ and $v$ denote the unknowns in \eqref{eq:coupled-to-coupled}.

Regarding well-posedness, Bona, Saut and Chen \cite{Bona-Saut-Chen-M-II} established local well-posedness for the corresponding Cauchy problem in $H^{s}(\mathbb{R})\times H^{s}(\mathbb{R})$ for $s>\frac{3}{4}$. On the periodic domain, Yang and Zhang \cite{Xin-Yang-KDV-KDV} established analytic local well-posedness for the system \eqref{eq:coupled-to-coupled} for $s\geq -\frac{1}{2}$. These results provide a foundation for the low-regularity analysis of the periodic Coupled KdV system.

We now turn to the smoothing properties of periodic dispersive equations. In this context, smoothing refers to the gain of spatial regularity exhibited by the nonlinear correction to the corresponding linear evolution. For the periodic KdV equation, Erdoğan and Tzirakis \cite{Erdogan-Tzirakis-global-sol-per} established a fundamental smoothing result. In particular,  their analysis yields a smoothing gain of order
\begin{equation*}
a<\min \{ 1,2s+1\}.
\end{equation*}
Their approach combines the Fourier restriction norm method with normal-form reductions obtained through differentiation by parts. Compaan \cite{E-compaan} subsequently investigated the smoothing properties of the Majda--Biello system on the torus. More recently, Başakoğlu and Gürel \cite{Engin-Burak} studied the smoothing properties of the Hirota--Satsuma system on the torus, contributing further results to the study of nonlinear smoothing for coupled KdV-type systems.
An observation by Oh \cite{Oh-diaphantine-cond-in-WPS} concerning the resonance structure of the Majda--Biello system provides further motivation for the present study. Oh observed that, when the coupling parameter $\alpha$ is negative, the relevant resonance denominator has a favorable nonvanishing structure away from the trivial resonances and yields additional frequency decay. In particular, when $\alpha =-1$, the smoothing estimate can be established by a comparatively direct application of KdV-type arguments, leading to the expected range
\begin{equation*}
a<\min\{1,2s+1\}.
\end{equation*}

Together with the appearance of terms analogous to those arising in the periodic KdV equation after normal-form reduction, this observation may initially suggest that establishing smoothing for the Coupled KdV system \eqref{eq:coupled-to-coupled} should be relatively straightforward. However, our analysis reveals that this expectation is not borne out: the interactions among the nonlinear terms introduce additional difficulties that are absent from the corresponding KdV and Majda--Biello problems. More precisely, the normal-form reduction of the first equation generates two terms involving the interaction $(u,u,v)$, namely $R_{3}(u,u,v)_{k}$ and $R_{7}(u,u,v)_{k}$, whose resonance structure presents a significant obstruction to the smoothing estimates. In certain high-low frequency regimes, the cubic contributions to the associated phase function can nearly cancel, leaving the phase bounded rather than providing the frequency-dependent lower bound needed to gain additional regularity. Consequently, the phase does not yield the same frequency gain as in the favorable resonance regime of the Majda--Biello system when $\alpha =-1$. The first equation also contains two further terms involving the interaction $(u,u,v)$, namely $R_{4}(u,u,v)_{k}$ and $R_{12}(u,u,v)_{k}$. Although their resonance functions can vanish in certain frequency regions, additional factors in their denominators compensate for the resulting loss of frequency decay. For $R_{3}$ and $R_{7}$, however, the interplay between the degeneracy of the phase and the accompanying frequency factors leads to the more restrictive smoothing range
\begin{equation*}
a<\min\left\{ 1,\dfrac{2}{3}(2s+1)\right\}.
\end{equation*}
Thus, despite the apparent similarities between the Coupled KdV system and each of the periodic KdV and Majda--Biello equations, the coupled nonlinear structure of the present system creates substantial additional difficulties in establishing smoothing.

Our proof combines the normal-form method, implemented through differentiation by parts \cite{BabinIlinTiti2011}, with duality arguments and estimates in Bourgain spaces. After the normal-form reduction, we organize the resulting terms according to the four interaction types
\begin{equation*}
(f,f,f), \qquad (f,f,g), \qquad (f,g,g), \qquad (g,g,g),
\end{equation*}
where $f$ and $g$ denote the interaction profiles associated with $u$ and $v$, respectively. This decomposition enables us to identify the relevant resonant sets and exploit the nonvanishing frequency factors imposed by the normal-form reduction.
In estimating the resulting terms $R_{i}$ and $\widetilde{R}_{i}$, we use the phase function
\begin{equation*}
\Phi=\sum \limits_{j=1}^{4}\phi^{\alpha_{j}}(k_{j}), \qquad \phi^{\alpha}(k)=\alpha k^{3},\qquad \alpha_{j}\in \left \{\dfrac{1}{6},-\dfrac{1}{6}\right\}.
\end{equation*}
The phase function appears in the denominator of the estimates arising from the duality argument. The restrictions defining the corresponding resonant sets guarantee that the relevant products of pairwise frequency sums are nonzero for the interactions under consideration. This nonvanishing property enables us to take advantage of the phase factor $\langle \Phi \rangle$ in the denominator and thereby obtain the necessary frequency decay. Combining normal-form restrictions with suitable low-high frequency decompositions and Bourgain-space estimates, we establish the desired nonlinear smoothing result for the periodic Coupled KdV system. This result shows that the difference between the nonlinear solution and its corresponding linear evolution belongs to a Sobolev space of higher regularity than that of the initial data. With these considerations, we are now in a position to state our principal result.

\begin{theorem}\label{thm:main}
    Fix $s>-\frac{1}{2}$. Let $(u,v)$ be the solution to the coupled system \eqref{eq:coupled-to-coupled} corresponding to the initial data $(u_{0},v_{0})\in H^{s}_{0}(\mathbb{T}) \times H^{s}_{0}(\mathbb{T}) $. Then, for every $a<\min \{ 1,\frac{4s}{3}+\frac{2}{3}\}$, the nonlinear parts of the solution satisfy 
\begin{equation}
u-e^{-t\frac{\partial_{xxx}}{6}}u_{0} \in C_{t}^{0}H^{a+s}_{x},
\end{equation}
and 
\begin{equation}
v-e^{t\frac{\partial_{xxx}}{6}}v_{0} \in C_{t}^{0}H^{a+s}_{x}.
\end{equation}
    
\end{theorem}
\begin{remark}
    Throughout this paper, we consider initial data in $H_{0}^{s}(\mathbb{T})\times H_{0}^{s}(\mathbb{T})$,  motivated by the presence of the terms $kf_{k}f_{0}$ and $kf_{k}g_{0}$ arising in the proofs of Propositions~\ref{eq:u_k} and~\ref{eq:v_k}.
\end{remark} 
\subsection{Notation}
Throughout this paper, we employ the notation $A \lesssim B$ to indicate that there exists a positive constant $C$, independent of the quantities under consideration, such that $A \leq C B$. The constant $C$ may vary from line to line and may depend on fixed parameters, such as the dimension $d$. Likewise, we write $A \ll B$ whenever $A \leq C^{-1} B$ for some sufficiently large constant $C$. Finally, the notation $A \sim B$ signifies that $A \lesssim B$ and $B \lesssim A$.
\begin{itemize}
\item The temporal Fourier transform and its inverse are defined standardly as follows.
For any function $f(t)$ on $\mathbb{R}$,
\begin{equation}
(\mathcal{F}_t f)(\tau)
:=
\int_{\mathbb{R}} e^{-it\tau} f(t)\, dt,
\qquad \forall \tau \in \mathbb{R}.
\end{equation}
On the other hand, for any function $g(\tau)$ on $\mathbb{R}$,
\begin{equation}
(\mathcal{F}_t^{-1} g)(t)
:=
\frac{1}{2\pi}
\int_{\mathbb{R}} e^{it\tau} g(\tau)\, d\tau,
\qquad \forall t \in \mathbb{R}.
\end{equation}
\item The spatial Fourier transform for any function $f(x)$ on $\mathbb{T}$,
\begin{equation}
(\mathcal{F}_x f)(k)
:=
\int_0^{2\pi}
e^{-ikx} f(x)\, dx,
\qquad \forall k \in \Z.
\end{equation}
On the other hand, for any function $\hat{f}(k)$ on $\Z$,
\begin{equation*}
(\mathcal{F}_x^{-1} f)(k)
:=\sum_{k\in \Z}
 \widehat{f}(k)e^{ikx}.
\end{equation*}
Consequently, the norm for the Sobolev space $H^s(\mathbb{T})$ is defined as
\begin{equation}
\Vert f\Vert_{H^{s}(\mathbb{T})}
=
\Vert \langle k \rangle^{s} \widehat{f}(k)\Vert_{\ell^2(\Z)},
\end{equation}
where $\langle k\rangle := \sqrt{1+|k|^{2}}$.\\
$\mathcal{H}^{s}:=H_{m_{1}}^{s}(\mathbb{T})\times H_{m_{2}}^{s}(\mathbb{T})$ and $H_{m}^{s}(\mathbb{T})=\big \{f\in H^{s}(\mathbb{T}): \int_{0}^{2\pi }f(x)dx=m \big \}$.
In particular, $H_{0}^{s}(\mathbb{T})$ is the homogeneous subspace of $H^{s}(\mathbb{T})$.
\item The space-time Fourier transform and its inverse are defined as
$\mathcal{F}=\mathcal{F}_{x}\mathcal{F}_{t}$
and
$\mathcal{F}^{-1}=\mathcal{F}_{x}^{-1}\mathcal{F}_{t}^{-1}$.
More precisely, for any function $f(x,t)$ on $\mathbb{T}\times\mathbb{R}$,
\begin{equation}
\widehat{f}(k,\tau)=(\mathcal{F}f)(k,\tau)
:=
\int_{\mathbb{R}}
\int_{0}^{2\pi}
e^{-i(t\tau+kx)} f(x,t)\, dx\, dt,
\qquad
\forall (k,\tau)\in \Z \times\mathbb{R}.
\end{equation}
For any function $g(k,\tau)$ on $\Z\times\mathbb{R}$,
\begin{equation}
\check{g}(x,t)=(\mathcal{F}^{-1} g)(x,t)
:=
\frac{1}{4\pi^{2}}
\int \limits_{\R}\sum \limits_{k\in \Z}
e^{i(t\tau+kx)} g(k,\tau)\, d\tau, \, 
\qquad
\forall (x,t)\in\mathbb{T} \times\mathbb{R}.
\end{equation}
\item The Fourier restriction spaces $X_{\alpha}^{s,b}$ associated to the system (\ref{eq:coupled-to-coupled}) are defined via the norm 
\begin{equation*}
\Vert u \Vert _{X_{\alpha}^{s,b}} =\Vert u \Vert _{X_{\phi^{\alpha}(k)}^{s,b}}= \Vert \langle k \rangle^{s}\langle \tau -\phi^{\alpha}(k) \rangle^{b} \hat{u}(k,\tau) \Vert_{\ell^{2}_{k}L^{2}_{\tau}},
\end{equation*}
where $\phi^{\alpha}(k)=\alpha k^{3}$ for some $\alpha \in \R.$\\
The restricted Bourgain norm is also defined by
\begin{equation*}
\Vert u \Vert _{X_{\alpha,\delta}^{s,b}}=\Vert u \Vert _{X_{\phi^{\alpha}(k),\delta}^{s,b}}=\inf \Big \{\Vert \tilde{u} \Vert _{X_{\phi^{\alpha}(k)}^{s,b}}:u=\tilde{u} \text{ on }[-\delta, \delta]\Big \}.
\end{equation*}
Observe that the parameters $\alpha $ corresponding to the first and second equations in the system (\ref{eq:coupled-to-coupled}) are given by $\frac{1}{6}$ and $-\frac{1}{6}$, respectively.
\item Finally, we use the subscript  $k$ to denote the Fourier coefficient of $R_{i}$. Accordingly, we express $R_{i}(u,v,w)$ through its Fourier series as
\begin{equation*}
R_{i}(u,v,w)=\sum \limits_{k}R_{i}(u,v,w)_{k}e^{ikx}.
\end{equation*}
The notation $\sum \limits^{\ast}$ indicates that the summation is restricted to terms for which the denominator does not vanish.
\end{itemize}
\section{Preliminaries}
\label{sec:prelim}
This section is devoted to recalling the preliminary estimates and well-posedness results that provide the analytical framework for establishing the main smoothing theorem. We begin by introducing two fundamental estimates in Bourgain spaces, namely the standard Duhamel estimate and the periodic $L^{6}$ Strichartz estimate, both of which play an essential role in the subsequent analysis. We then state the local well-posedness results on both $\R$ and $\T$. 
\begin{lemma}\cite{Ginibre-Velo}\label{eq:generic-lemma}
    For $\frac{1}{2}<b\leq 1$ and $\alpha \in \R$ with $\alpha\neq 0$, we have 
     \begin{equation*}
\Bigg \Vert\psi_{\delta}(t)\int \limits_{0}^{t}e^{(t'-t)(\alpha\partial_{xxx})}F(t')dt' \Bigg \Vert_{X_{\tau=\phi^{\alpha}(k)}^{s,b}}\lesssim \Vert F \Vert_{X_{\tau=\phi^{\alpha}(k), \, \delta}^{s,b-1}} ,
\end{equation*}
  where $\phi^{\alpha}(k)=\alpha k^{3}$.\end{lemma}
\begin{proposition}\label{eq:inverse-fourier-L6-norm}\cite{Bourgain-partII}
    For any $\epsilon>0$ and $b>\frac{1}{2}$, we have
    \begin{equation*}
\left\| \chi_{[-\delta,\delta]}(t)\, u \right\|_{L^6_{t,x}(\mathbb{R}\times \mathbb{T})}
\leq C_{\varepsilon,b}\, \|u\|_{X^{\varepsilon,b}_\delta}.
\end{equation*}
\end{proposition}
We next recall the well-posedness results relevant to the coupled KdV system under consideration.
Bona, Chen, and Saut \cite{Bona-Saut-Chen-M-II} established local well-posedness for the following rescaled coupled KdV system posed on the real line:
\begin{equation}
\label{eq:bona-chen-saut}
\begin{cases}
        u_t+v_x+(uv)_x+v_{xxx}=0,\\[2mm]
        v_t+u_x+vv_x+u_{xxx}=0,\\[2mm]
        (u,v)|_{t=0}=(u_0,v_0), 
    \end{cases}
\end{equation}
where $(x,t)\in\mathbb{R}\times\mathbb{R}$. Their result is stated as follows.
\begin{theorem}\cite{Bona-Saut-Chen-M-II}
    Let $u_{0},v_{0}\in H^{s}(\mathbb{R})$ with $s>\frac{3}{4}$. Then there exists a value $T>0$ such that the system (\ref{eq:bona-chen-saut}) has a unique solution $(u,v)$ in $C(0,T;H^{s}(\mathbb{R}))^{2}$.
\end{theorem}
\noindent In the periodic setting, Yang and Zhang \cite{Xin-Yang-KDV-KDV} established the following local well-posedness result.
\begin{theorem}\cite{Xin-Yang-KDV-KDV}
 Let   $\mathcal{H}^{s}(\mathbb{T})=H_{m_{1}}^{s}(\mathbb{T})\times H_{m_{2}}^{s}(\mathbb{T}) $. Then (\ref{eq:coupled-to-coupled}) is analytically locally well posed in $\mathcal{H}^{s}(\mathbb{T})$ for $s\geq -\frac{1}{2}$.
\end{theorem}
\section{Proof of Theorem \ref{thm:main}}
\label{sec:strips}
We begin by presenting the principal steps in the proof of Theorem \ref{thm:main}. The proof relies on the normal-form method, combined with multilinear estimates in Bourgain spaces and a detailed analysis of the associated frequency interactions. We first apply differentiation by parts in time to the nonlinear terms of the system \eqref{eq:coupled-to-coupled}. This procedure allows us to express the nonlinear contributions in terms of  bilinear terms and higher-order multilinear expressions. We then decompose each multilinear term into resonant and nonresonant contributions, thereby exploiting the structure of the associated resonance functions appearing in the denominator. 

The resonant and bilinear terms are estimated directly in the $L^{\infty }_{t}H^{s+a}_{x}$ norm, as established in Propositions \ref{prop:prop-B-1-2-3}-\ref{prop:mathcal-R-tilde}, whereas the resulting nonresonant terms are controlled through duality arguments in Bourgain spaces and suitable frequency decompositions, as demonstrated in Propositions \ref{eq:burak-tzirakis-R1}-\ref{prop:R-tilde-5-8}. A central role in this analysis is played by the resonance function $\Phi$
which arises in the frequency-space formulation of the multilinear estimates. The algebraic structure of this function is essential for obtaining the lower bounds required to establish the desired smoothing estimates.
This phase function arises in the denominator after applying the duality argument.  To exploit the nonvanishing properties of the interactions, we decompose the corresponding contributions into  four interaction types
 \begin{equation*}
(f,f,f),\quad (f,f,g), \quad (f,g,g), \quad (g,g,g).
\end{equation*} 
 The restrictions imposed by $\mathcal{R}_{f,f,f}$, $\mathcal{R}_{f,f,g}$, $\mathcal{R}_{g,g,g}$ ensure that \begin{equation*}
(k_{1}+k_{2})(k_{1}+k_{3})(k_{2}+k_{3})\neq 0
\end{equation*} for the terms $R_{i}$, where $i\in\{1,3,4,5,7,8,10,11,12\}$. Likewise, the restriction imposed by $\mathcal{R}_{f,g,g}$ guarantees that \begin{equation*}
(k_{1}+k_{3})(k_{2}+k_{3})\neq 0
\end{equation*} for the terms $R_{6}$ and $R_{9}$. The argument used in Proposition \ref{prop:R-2} applies directly to all $(f,g,g)$-type terms satisfying the non-resonance condition $k_{2}+k_{3}\neq 0$. Therefore, it remains only to analyze the $R_{6}$ and $R_{9}$ terms for which the resonance condition $k_{2}+k_{3}=0$ may occur, in view of the factorization of the phase function
\begin{equation*}
\left \vert \Phi \right \vert \sim \left \vert (k_{2}+k_{3})\big[2(k_{2}+k_{3})^{2}-3(k_{1}k_{4}+k_{2}k_{3})\big]\right \vert.
\end{equation*} These non-vanishing conditions will play a crucial role in the subsequent propositions, as they allow us to exploit the presence of the factor $\langle  \Phi  \rangle$ in the denominator. 
We next turn to the analysis of the resonant regions associated with the $\widetilde{R}_{i}$-terms. These regions are characterized by the conditions
\begin{equation*}
(k_{1}+k_{2})(k_{1}+k_{3})( k_{2}+k_{3})  =0 \text{ and }(k_{1}+k_{2})( k_{2}+k_{3})  =0
\end{equation*} 
 corresponding to the terms $\widetilde{R}_{1},\widetilde{R}_{2},\widetilde{R}_{4},\widetilde{R}_{5},\widetilde{R}_{6},\widetilde{R}_{7},\widetilde{R}_{8},\widetilde{R}_{9},\widetilde{R}_{10},\widetilde{R}_{11}$, and $\widetilde{R}_{3}$, $\widetilde{R}_{12}$, respectively. The analysis follows the same approach as in the preceding argument.

 As noted in the proof of Proposition \ref{eq:u_k}, we assume throughout the remainder of the paper that the frequencies $k$, $k_{1}$, $k_{2}$, $k_{3}$, and $k_{4}$ are all nonzero.

\subsection{Normal form reduction}
\label{sec:normal-form}
In this subsection, we apply integration by parts in time to each equation of the system \eqref{eq:coupled-to-coupled} to derive the corresponding normal-form representations. To take advantage of the phase functions appearing in the denominators of the resulting multilinear terms, we decompose these terms into resonant contributions and remainder terms $R_{i}-\wtr_{i}$. 
\begin{proposition}\label{eq:u_k}
    Let $(u_{0},v_{0}) \in H^{s}_{0}(\mathbb{T})\times  H^{s}_{0}(\mathbb{T})$. Then we have 
    \begin{align}
\partial_t \Big[
e^{-it\frac{k^{3}}{6}}
\big(
u_k
- B_1(u,u)_k
+ B_2(v,v)_k
-& B_3(u,v)_k
\big)
\Big]\notag
=
e^{-it\frac{k^{3}}{6}}
\Big[
R_1(u,u,u)_k
+ R_2(u,v,v)_k\notag\\
&\qquad + R_3(u,u,v)_k+ R_4(u,u,v)_k
+ R_5(v,v,v)_k\notag\\
&\qquad + R_6(u,v,v)_k+ R_{7}(u,u,v)_k  
+ R_{8}(v,v,v)_k\notag\\
&\qquad + R_{9}(u,v,v)_k+ R_{10}(u,u,u)_k  
+ R_{11}(u,v,v)_k\notag\\
&\qquad + R_{12}(u,u,v)_k+\mathcal{R}_{(u,u,u)_{k}}+\mathcal{R}_{(u,u,v)_{k}}+\mathcal{R}_{(v,v,v)_{k}}\notag\\
& \qquad +\mathcal{R}_{(u,v,v)_{k}}\Big],
\label{eq:normal-first-eq}
\end{align}
where
\begin{gather*}
B_1(u,u)_k=\sum \limits_{k_{1}+k_{2}=k}^{\ast}\frac{9k}{2}\frac{u_{k_{1}}u_{k_{2}}}{k^{3}-k_{1}^{3}-k_{2}^{3}},\\
B_2(v,v)_k=\sum \limits_{k_{1}+k_{2}=k}^{\ast}\frac{3k}{2}\frac{v_{k_{1}}v_{k_{2}}}{k^{3}+k_{1}^{3}+k_{2}^{3}},\\
B_3(u,v)_k= \sum \limits_{k_{1}+k_{2}=k}^{\ast}3k\frac{u_{k_{1}}v_{k_{2}}}{k^{3}-k_{1}^{3}+k_{2}^{3}},\\
R_1(u,u,u)_k=\frac{27i}{4}\sum \limits_{\substack{ k_{1}+k_{2}+k_{3}=k    \\   (k_{1}+k_{2})(k_{1}+k_{3})( k_{2}+k_{3})  \neq 0    }}^{\ast}\frac{ k(k_{2}+k_{3})\, u_{k_{1}}u_{k_{2}}u_{k_{3}}}{k^{3}-k_{1}^{3}-(k_{2}+k_{3})^{3}},\\
R_2(u,v,v)_k=-\frac{9i}{4}\sum \limits_{\substack{ k_{1}+k_{2}+k_{3}=k    \\   k_{2}+k_{3}  \neq 0 }}^{\ast}\frac{ k(k_{2}+k_{3})\, u_{k_{1}}v_{k_{2}}v_{k_{3}}}{k^{3}-k_{1}^{3}-(k_{2}+k_{3})^{3}},\\
R_3(u,u,v)_k=\frac{9i}{2}\sum \limits_{\substack{ k_{1}+k_{2}+k_{3}=k    \\   (k_{1}+k_{2})(k_{1}+k_{3})( k_{2}+k_{3})  \neq 0 }}^{\ast}\frac{k(k_{2}+k_{3}) \, u_{k_{1}}u_{k_{2}}v_{k_{3}}}{k^{3}-k_{1}^{3}-(k_{2}+k_{3})^{3}},\\
R_4(u,u,v)_k=\frac{3i}{4}\sum \limits_{\substack{ k_{1}+k_{2}+k_{3}=k    \\   (k_{1}+k_{2})(k_{1}+k_{3})( k_{2}+k_{3})  \neq 0}}^{\ast}\frac{k(k_{1}+k_{2})\, u_{k_{1}}u_{k_{2}}v_{k_{3}}}{k^{3}+(k_{1}+k_{2})^{3}+k_{3}^{3}},\\
R_5(v,v,v)_k=-\frac{9i}{4}\sum \limits_{\substack{  k_{1}+k_{2}+k_{3}=k   \\    (k_{1}+k_{2})(k_{1}+k_{3})(k_{2}+k_{3})  \neq 0 }}^{\ast}\frac{k(k_{2}+k_{3})\, v_{k_{1}}v_{k_{2}}v_{k_{3}}}{k^{3}+k_{1}^{3}+(k_{2}+k_{3})^{3}},\\
R_6(u,v,v)_k=-\frac{3i}{2}\sum \limits_{\substack{  k_{1}+k_{2}+k_{3}=k   \\   (k_{1}+k_{3})(k_{2}+k_{3})  \neq 0 }}^{\ast}\frac{k(k_{1}+k_{3})\, u_{k_{1}}v_{k_{2}}v_{k_{3}}}{k^{3}+(k_{1}+k_{3})^{3}+k_{2}^{3}},\\
R_{7}(u,u,v)_k=\frac{9i}{4} \sum \limits_{\substack{ k_{1}+k_{2}+k_{3}=k    \\   (k_{1}+k_{2})(k_{1}+k_{3})( k_{2}+k_{3})  \neq 0}}^{\ast}\frac{k(k_{1}+k_{2})\, u_{k_{1}}u_{k_{2}}v_{k_{3}}}{k^{3}-(k_{1}+k_{2})^{3}+k_{3}^{3}},\\
R_{8}(v,v,v)_k=-\frac{3i}{4} \! \sum \limits_{\substack{ k_{1}+k_{2}+k_{3}=k    \\   (k_{1}+k_{2})(k_{1}+k_{3}) (k_{2}+k_{3})  \neq 0 }}^{\ast}\frac{k(k_{2}+k_{3})\, v_{k_{1}}v_{k_{2}}v_{k_{3}}}{k^{3}+k_{1}^{3}-(k_{2}+k_{3})^{3}},\\
R_{9}(u,v,v)_k=\frac{3i}{2}  \sum \limits_{\substack{ k_{1}+k_{2}+k_{3}=k    \\   (k_{1}+k_{3})(k_{2}+k_{3})  \neq 0 }}^{\ast}\frac{k(k_{1}+k_{3})\, u_{k_{1}}v_{k_{2}}v_{k_{3}}}{k^{3}-(k_{1}+k_{3})^{3}+k_{2}^{3}},\\
R_{10}(u,u,u)_k=-\frac{3i}{4} \! \sum \limits_{\substack{  k_{1}+k_{2}+k_{3}=k   \\   (k_{1}+k_{2})(k_{1}+k_{3})( k_{2}+k_{3})  \neq 0 }}^{\ast}\frac{k(k_{2}+k_{3})\, u_{k_{1}}u_{k_{2}}u_{k_{3}}}{k^{3}-k_{1}^{3}+(k_{2}+k_{3})^{3}},\\
R_{11}(u,v,v)_k=\frac{9i}{4} \! \sum \limits_{\substack{    k_{1}+k_{2}+k_{3}=k \\   k_{2}+k_{3}  \neq 0 }}^{\ast}\frac{k(k_{2}+k_{3})\, u_{k_{1}}v_{k_{2}}v_{k_{3}}}{k^{3}-k_{1}^{3}+(k_{2}+k_{3})^{3}},\\
R_{12}(u,u,v)_k=\frac{3i}{2}  \sum \limits_{\substack{   k_{1}+k_{2}+k_{3}=k  \\  (k_{1}+k_{2})(k_{1}+k_{3})( k_{2}+k_{3})  \neq 0  }}^{\ast}\frac{k(k_{2}+k_{3})\, u_{k_{1}}u_{k_{2}}v_{k_{3}}}{k^{3}-k_{1}^{3}+(k_{2}+k_{3})^{3}},\\
\mathcal{R}_{(u,u,u)_{k}}=-\frac{12i}{5}\frac{\vert u_{k}\vert^{2}u_{k}}{k}-\frac{3i}{2}ku_{k}\sum \limits_{\vert j \vert \neq \vert k \vert}\frac{\vert u_{j}\vert^{2}}{2k^{2}-kj+2j^{2}},
\end{gather*}
\begin{align*}
\mathcal{R}_{(u,u,v)_{k}} &={} -\frac{6i}{5} \frac{\vert u_{k}\vert^{2}v_{k}}{k} -\frac{3i}{8} \frac{ u_{k} u_{k}v_{-k}}{k} +\frac{3i}{2}u_{k}\sum \limits_{\vert j \vert \neq \vert k\vert}^{\ast}\frac{u_{j}\overline{v_{j}}}{j} +\frac{3i}{2}u_{k}\sum \limits_{\vert j \vert \neq \vert k\vert}^{\ast}\frac{(j+k)u_{j}\overline{v_{j}}}{2k^{2}+3jk+3j^{2}}\\
&\quad{}-\frac{9i}{2}ku_{k}\sum \limits_{\vert j \vert \neq \vert k\vert}^{\ast}\frac{(j+k)u_{j}\overline{v_{j}}}{j(2j^{2}+3kj+3k^{2})} +\frac{3i}{2}kv_{k}\sum \limits_{\vert j \vert \neq \vert k\vert}^{\ast}\frac{\vert u_{j}\vert^{2}}{2k^{2}-kj+2j^{2}}\\
&\quad{}+\frac{3i}{2}ku_{k}\sum \limits_{\vert j \vert \neq \vert k\vert}^{\ast}\frac{ u_{j}\overline{v_{j}}}{2k^{2}-kj+2j^{2}},
\end{align*}
\begin{align*}
\mathcal{R}_{(v,v,v)_{k}} &={} -\frac{3i}{8}\frac{\vert v_{k}\vert^{2}v_{k}}{k} -\frac{9i}{2}v_{k}\sum \limits_{\vert j \vert \neq  \vert k \vert}\frac{(k-j)\vert v_{j}\vert^{2}}{2k^{2}-3kj+3j^{2}}\\
&\quad{}+6ikv_{k}\sum \limits_{\substack{ j\geq 1\\ j\neq \vert k\vert}}^{\ast}\frac{j^{2}\vert v_{j}\vert^{2}}{(2j^{2}-3kj+3k^{2})(2j^{2}+3kj+3k^{2})},
\end{align*}
\begin{align*}
\mathcal{R}_{(u,v,v)_{k}} &={} -\frac{3i}{4}\frac{u_{k}\vert v_{k}\vert^{2}}{k} -\frac{3i}{2}u_{k}\sum \limits_{\vert j \vert \neq  \vert k \vert}^{\ast}\frac{(k-j)\vert v_{j}\vert^{2}}{2k^{2}-3kj+3j^{2}}\\
&\quad{}-6iku_{k}\sum \limits_{ \substack{j\geq 1  \\ j\neq \vert k \vert} }^{\ast}\frac{j^{2}\vert v_{j}\vert^{2}}{(2j^{2}-3kj+3k^{2})(2j^{2}+3kj+3k^{2})}.
\end{align*}
\end{proposition}
\begin{proof}
    Taking the Fourier transform of (\ref{eq:coupled-to-coupled}) with respect to the space variable, we obtain the following equations
    \begin{gather*}
\partial_{t}u_{k}-\frac{ik^{3}}{6}u_{k}+\frac{3i}{4}k\sum \limits_{k_{1}+k_{2}=k}u_{k_{1}}u_{k_{2}}-\frac{i}{4}k\sum \limits_{k_{1}+k_{2}=k}v_{k_{1}}v_{k_{2}}+\frac{i}{2}k\sum \limits_{k_{1}+k_{2}=k}u_{k_{1}}v_{k_{2}}=0,\\
\partial_{t}v_{k}+\frac{ik^{3}}{6}v_{k}-\frac{i}{4}k\sum \limits_{k_{1}+k_{2}=k}u_{k_{1}}u_{k_{2}}+\frac{3i}{4}k\sum \limits_{k_{1}+k_{2}=k}v_{k_{1}}v_{k_{2}}+\frac{i}{2}k\sum \limits_{k_{1}+k_{2}=k}u_{k_{1}}v_{k_{2}}=0.
\end{gather*}
We now introduce the interaction profiles $f_{k}(t)=e^{-it\frac{k^{3}}{6}}u_{k}(t)$ and $g_{k}(t)=e^{it\frac{k^{3}}{6}}v_{k}(t)$ to the above equations.
\begin{align*}
e^{\frac{itk^{3}}{6}}\partial_{t}f_{k} &={} -\frac{3i}{4}k\sum \limits_{k_{1}+k_{2}=k}e^{-it\frac{-k_{1}^{3}-k_{2}^{3}}{6}}f_{k_{1}}f_{k_{2}} +\frac{i}{4}k\sum \limits_{k_{1}+k_{2}=k}e^{-it\frac{k_{1}^{3}+k_{2}^{3}}{6}}g_{k_{1}}g_{k_{2}}\\
&\quad{}-\frac{i}{2}k\sum \limits_{k_{1}+k_{2}=k}e^{-it\frac{-k_{1}^{3}+k_{2}^{3}}{6}}f_{k_{1}}g_{k_{2}},
\end{align*}
\begin{align*}
e^{-\frac{itk^{3}}{6}}\partial_{t}g_{k} &={} \frac{i}{4}k\sum \limits_{k_{1}+k_{2}=k}e^{-it\frac{-k_{1}^{3}-k_{2}^{3}}{6}}f_{k_{1}}f_{k_{2}} -\frac{3i}{4}k\sum \limits_{k_{1}+k_{2}=k}e^{-it\frac{k_{1}^{3}+k_{2}^{3}}{6}}g_{k_{1}}g_{k_{2}}\\
&\quad{}-\frac{i}{2}k\sum \limits_{k_{1}+k_{2}=k}e^{-it\frac{-k_{1}^{3}+k_{2}^{3}}{6}}f_{k_{1}}g_{k_{2}}.
\end{align*}
So, we have
\begin{align}
\partial_{t}f_{k} &={} -\frac{3i}{4}k\sum \limits_{k_{1}+k_{2}=k}e^{-it\frac{k^{3}-k_{1}^{3}-k_{2}^{3}}{6}}f_{k_{1}}f_{k_{2}} +\frac{i}{4}k\sum \limits_{k_{1}+k_{2}=k}e^{-it\frac{k^{3}+k_{1}^{3}+k_{2}^{3}}{6}}g_{k_{1}}g_{k_{2}}\notag\\
&\quad{}-\frac{i}{2}k\sum \limits_{k_{1}+k_{2}=k}e^{-it\frac{k^{3}-k_{1}^{3}+k_{2}^{3}}{6}}f_{k_{1}}g_{k_{2}},
\label{eq:f_k}
\end{align}
\begin{align}
\partial_{t}g_{k} &={} \frac{i}{4}k\sum \limits_{k_{1}+k_{2}=k}e^{-it\frac{-k^{3}-k_{1}^{3}-k_{2}^{3}}{6}}f_{k_{1}}f_{k_{2}} -\frac{3i}{4}k\sum \limits_{k_{1}+k_{2}=k}e^{-it\frac{-k^{3}+k_{1}^{3}+k_{2}^{3}}{6}}g_{k_{1}}g_{k_{2}}\notag\\
&\quad{}-\frac{i}{2}k\sum \limits_{k_{1}+k_{2}=k}e^{-it\frac{-k^{3}-k_{1}^{3}+k_{2}^{3}}{6}}f_{k_{1}}g_{k_{2}}.
\label{eq:g_k}
\end{align}

We begin by analyzing the resonant contributions. Since all terms vanish when $k=0$, it suffices to restrict our attention to $k\neq 0$.
For the first term in (\ref{eq:f_k}), the resonance condition is $kk_{1}k_{2}=0$. Since $k\neq 0$, we distinguish the following three regions:
\begin{align*}
S_{1}&=\{k_{1}=0,\ k_{2}=0\},\\
S_{2}&=\{k_{1}=0,\ k_{2}\neq 0\},\\
S_{3}&=\{k_{1}\neq 0,\ k_{2}=0\}.
\end{align*}
The region $S_{1}$ is excluded, since $k_{1}=k_{2}=0$ implies $k=0$, contradicting our assumption. In $S_{2}$, the condition $k_{1}=0$ yields $k_{2}=k$, whereas in $S_{3}$, the condition $k_{2}=0$ gives $k_{1}=k$. Consequently, both regions contribute equally, namely,
\begin{equation*}
\left.\text{First term}\right|_{S_{2}}
=\left.\text{First term}\right|_{S_{3}}
=-\frac{3i}{4}kf_{k}f_{0}.
\end{equation*}
Combining these contributions, the resonant part of the first term is
\begin{equation*}
-\frac{3i}{2}kf_{k}f_{0}.
\end{equation*}
For the second term in (\ref{eq:f_k}), the resonance condition takes the form
$k(3k_{1}^{2}-3kk_{1}+2k^{2})=0.$
Since $k\neq 0$ and the quadratic polynomial $3k_{1}^{2}-3kk_{1}+2k^{2}$ has a negative discriminant, this condition admits no real solution. Thus, the second term has no resonant contribution.
Finally, for the third term in (\ref{eq:f_k}), the resonance condition is
$k_{2}(2k_{2}^{2}-3kk_{2}+3k^{2})=0.$
The quadratic factor is strictly positive for $k\neq 0$, so resonance occurs only when $k_{2}=0$. This yields the contribution
\begin{equation*}
-\frac{i}{2}kf_{k}g_{0}.
\end{equation*}
Therefore, collecting the resonant contributions from all three terms, we obtain the residual term associated with (\ref{eq:f_k}):
\begin{equation*}
-\frac{3i}{2}kf_{k}f_{0}
-\frac{i}{2}kf_{k}g_{0}.
\end{equation*}
Since the terms $kf_{k}f_{0}$ and $kf_{k}g_{0}$ cannot be controlled within the framework of our method, we restrict the analysis to mean-zero functions, that is, $u_{0}=v_{0}=0$. Hence, there is no resonant term corresponding to the equation $\partial_{t}f_{k}$.
We now apply differentiation by parts to (\ref{eq:f_k}).
For the first term of (\ref{eq:f_k}), we utilize the identity 
\begin{equation*}
\partial_{t}\Bigg( \frac{e^{-it \frac{k^{3}-k_{1}^{3}-k_{2}^{3}}{6}  }}{-i\frac{k^{3}-k_{1}^{3}-k_{2}^{3}}{6}   } \Bigg)=e^{-it\frac{k^{3}-k_{1}^{3}-k_{2}^{3}}{6}}
\end{equation*}
to get
\begin{align*}
    -\frac{3i}{4}k\sum \limits_{k_{1}+k_{2}=k}^{\ast}&e^{-it\frac{k^{3}-k_{1}^{3}-k_{2}^{3}}{6}}f_{k_{1}}f_{k_{2}}\\
    &=\partial_{t}\Bigg(\sum \limits_{k_{1}+k_{2}=k}^{\ast} \frac{9}{2}k \frac{e^{-it\frac{k^{3}-k_{1}^{3}-k_{2}^{3}}{6}   }}{  k^{3}-k_{1}^{3}-k_{2}^{3}}f_{k_{1}}f_{k_{2}}\Bigg)-\sum \limits_{k_{1}+k_{2}=k}^{\ast}\frac{9}{2}k \frac{e^{-it\frac{k^{3}-k_{1}^{3}-k_{2}^{3}}{6}   }}{  k^{3}-k_{1}^{3}-k_{2}^{3} }\partial_{t}\Big( f_{k_{1}}f_{k_{2}}\Big)\\
    &:=\partial_{t}\Bigg( \sum \limits_{k_{1}+k_{2}=k}^{\ast} \frac{9}{2}k \frac{e^{-it\frac{k^{3}-k_{1}^{3}-k_{2}^{3}}{6}   }}{ k^{3}-k_{1}^{3}-k_{2}^{3}}f_{k_{1}}f_{k_{2}}\Bigg)+(I),
\end{align*}
where
\begin{align*}
(I)&=-\frac{9}{2} \sum \limits_{k_{1}+k_{2}=k}^{\ast}k \frac{e^{-it\frac{k^{3}-k_{1}^{3}-k_{2}^{3}}{6}   }}{ k^{3}-k_{1}^{3}-k_{2}^{3}}\Bigg[ \Bigg( -\frac{3i}{4}k_{1}\sum \limits_{k_{3}+k_{4}=k_{1}}^{\ast}e^{-it\frac{k_{1}^{3}-k_{3}^{3}-k_{4}^{3}}{6}}f_{k_{3}}f_{k_{4}}\\
&+\frac{i}{4}k_{1}\sum \limits_{k_{3}+k_{4}=k_{1}}^{\ast}e^{-it\frac{k_{1}^{3}+k_{3}^{3}+k_{4}^{3}}{6}}g_{k_{3}}g_{k_{4}}-\frac{i}{2}k_{1}\sum \limits_{k_{3}+k_{4}=k_{1}}^{\ast}e^{-it\frac{k_{1}^{3}-k_{3}^{3}+k_{4}^{3}}{6}}f_{k_{3}}g_{k_{4}}       \Bigg)f_{k_{2}}\\
&+\Bigg( -\frac{3i}{4}k_{2}\sum \limits_{k_{3}+k_{4}=k_{2}}^{\ast}e^{-it\frac{k_{2}^{3}-k_{3}^{3}-k_{4}^{3}}{6}}f_{k_{3}}f_{k_{4}}+\frac{i}{4}k_{2}\sum \limits_{k_{3}+k_{4}=k_{2}}^{\ast}e^{-it\frac{k_{2}^{3}+k_{3}^{3}+k_{4}^{3}}{6}}g_{k_{3}}g_{k_{4}}\\
&-\frac{i}{2}k_{2}\sum \limits_{k_{3}+k_{4}=k_{2}}^{\ast}e^{-it\frac{k_{2}^{3}-k_{3}^{3}+k_{4}^{3}}{6}}f_{k_{3}}g_{k_{4}}       \Bigg)f_{k_{1}}
         \Bigg]\\
&=\frac{27i}{8} \sum \limits_{\substack{ k_{2}+k_{3}+k_{4}=k    \\  k_{3}+k_{4}\neq 0    }}^{\ast} k(k_{3}+k_{4})\frac{e^{-it \frac{k^{3}-k_{2}^{3}-k_{3}^{3}-k_{4}^{3}}{6}   }}{   k^{3}-k_{2}^{3}-(k_{3}+k_{4})^{3}}f_{k_{2}}f_{k_{3}}f_{k_{4}}\\
&- \frac{9i}{8}\sum \limits_{\substack{  k_{2}+k_{3}+k_{4}=k   \\   k_{3}+k_{4}  \neq 0 }}^{\ast} k(k_{3}+k_{4})\frac{e^{-it\frac{k^{3}-k_{2}^{3}+k_{3}^{3}+k_{4}^{3}}{6}  }}{   k^{3}-k_{2}^{3}-(k_{3}+k_{4})^{3}}f_{k_{2}}g_{k_{3}}g_{k_{4}}\\
&+\frac{9i}{4}  \sum \limits_{\substack{  k_{2}+k_{3}+k_{4}=k    \\   k_{3}+k_{4}  \neq 0 }}^{\ast} k(k_{3}+k_{4})\frac{e^{-it\frac{k^{3}-k_{2}^{3}-k_{3}^{3}+k_{4}^{3}}{6}   }}{ k^{3}-k_{2}^{3}-(k_{3}+k_{4})^{3}}f_{k_{2}}f_{k_{3}}g_{k_{4}}\\
&+\frac{27i}{8} \sum \limits_{ \substack{   k_{1}+k_{3}+k_{4}=k   \\   k_{3}+k_{4}  \neq 0 }}^{\ast}  k(k_{3} + k_{4})\frac{e^{-it\frac{k^{3}-k_{1}^{3}-k_{3}^{3}-k_{4}^{3}}{6}   }}{ k^{3}-k_{1}^{3}-(k_{3}+k_{4})^{3} }f_{k_{1}}f_{k_{3}}f_{k_{4}}\\
&-\frac{9i}{8} \sum \limits_{ \substack{  k_{1}+k_{3}+k_{4}=k   \\   k_{3}+k_{4}  \neq 0 } }^{\ast}  k(k_{3} + k_{4})\frac{e^{-it\frac{k^{3}-k_{1}^{3}+k_{3}^{3}+k_{4}^{3}}{6}   }}{  k^{3}-k_{1}^{3}-(k_{3}+k_{4})^{3}}f_{k_{1}}g_{k_{3}}g_{k_{4}}\\
&+\frac{9i}{4} \sum \limits_{\substack{ k_{1}+k_{3}+k_{4}=k    \\   k_{3}+k_{4}  \neq 0 }  }^{\ast}  k(k_{3} + k_{4})\frac{e^{-it\frac{k^{3}-k_{1}^{3}-k_{3}^{3}+k_{4}^{3}}{6}   }}{  k^{3}-k_{1}^{3}-(k_{3}+k_{4})^{3} }f_{k_{1}}f_{k_{3}}g_{k_{4}}.
\end{align*}
After a suitable relabeling of the frequency variables, the first term in the expression $\partial_{t}f_{k}$ takes the form 
\begin{align}
-\frac{3i}{4}k\sum \limits_{k_{1}+k_{2}=k}^{\ast}e^{-it\Big(-k+k_{1}+k_{2}+\frac{k^{3}-k_{1}^{3}-k_{2}^{3}}{6}\Big)}&f_{k_{1}}f_{k_{2}}\notag
 =\partial_{t}\Bigg( \sum \limits_{k_{1}+k_{2}=k}^{\ast} \frac{9}{2}k \frac{e^{-it\frac{k^{3}-k_{1}^{3}-k_{2}^{3}}{6}   }}{  k^{3}-k_{1}^{3}-k_{2}^{3} }f_{k_{1}}f_{k_{2}}\Bigg)\notag\\
& +\frac{27i}{4} \sum \limits_{\substack{  k_{1}+k_{2}+k_{3}=k   \\   k_{2}+k_{3}  \neq 0 }}^{\ast} k(k_{2}+k_{3})\frac{e^{-it\frac{k^{3}-k_{1}^{3}-k_{2}^{3}-k_{3}^{3}}{6}  }}{  k^{3}-k_{1}^{3}-(k_{2}+k_{3})^{3} }f_{k_{1}}f_{k_{2}}f_{k_{3}}\notag\\
& - \frac{9i}{4}\sum \limits_{\substack{  k_{1}+k_{2}+k_{3}=k   \\   k_{2}+k_{3}  \neq 0 }}^{\ast} k(k_{2}+k_{3})\frac{e^{-it\frac{k^{3}-k_{1}^{3}+k_{2}^{3}+k_{3}^{3}}{6}  }}{  k^{3}-k_{1}^{3}-(k_{2}+k_{3})^{3} }f_{k_{1}}g_{k_{2}}g_{k_{3}}\notag\\
& +\frac{9i}{2}  \sum \limits_{\substack{   k_{1}+k_{2}+k_{3}=k   \\   k_{2}+k_{3}  \neq 0 }}^{\ast} k(k_{2} + k_{3})\frac{e^{-it\frac{k^{3}-k_{1}^{3}-k_{2}^{3}+k_{3}^{3}}{6}  }}{  k^{3}-k_{1}^{3}-(k_{2}+k_{3})^{3} }f_{k_{1}}f_{k_{2}}g_{k_{3}}.
\label{eq:fk-first}
\end{align}
For the second  and third terms of (\ref{eq:f_k}), we perform a differentiation by parts in time, following the same procedure as in the preceding computation. Consequently, the second and third terms in $\partial_{t}f_{k}$ admit the following representations:
\begin{align}
\frac{i}{4}k\sum \limits_{k_{1}+k_{2}=k}^{\ast}e^{-it\Big(-k-k_{1}-k_{2}+\frac{k^{3}+k_{1}^{3}+k_{2}^{3}}{6}\Big)}&g_{k_{1}}g_{k_{2}}
 =-\partial_{t}\Bigg(\sum \limits_{k_{1}+k_{2}=k}^{\ast} \frac{3k}{2} \frac{e^{-it\frac{k^{3}+k_{1}^{3}+k_{2}^{3}}{6}   }}{ k^{3}+k_{1}^{3}+k_{2}^{3}  }g_{k_{1}}g_{k_{2}}\Bigg)\notag\\
& +\frac{3i}{4} \sum \limits_{\substack{  k_{1}+k_{2}+k_{3}=k   \\   k_{1}+k_{2}  \neq 0 }}^{\ast} k(k_{1}+k_{2})\frac{e^{-it\frac{k^{3}-k_{1}^{3}-k_{2}^{3}+k_{3}^{3}}{6}  }}{ k^{3}+(k_{1}+k_{2})^{3}+k_{3}^{3} }f_{k_{1}}f_{k_{2}}g_{k_{3}}\notag\\
&- \frac{9i}{4}\sum \limits_{\substack{  k_{1}+k_{2}+k_{3}=k   \\   k_{2}+k_{3}  \neq 0 }}^{\ast} k(k_{2}+k_{3})\frac{e^{-it\frac{k^{3}+k_{1}^{3}+k_{2}^{3}+k_{3}^{3}}{6}   }}{  k^{3}+k_{1}^{3}+(k_{2}+k_{3})^{3} }g_{k_{1}}g_{k_{2}}g_{k_{3}}\notag\\
& -\frac{3i}{2}  \sum \limits_{\substack{  k_{1}+k_{2}+k_{3}=k   \\   k_{1}+k_{3}  \neq 0 } }^{\ast}k(k_{1} + k_{3})\frac{e^{-it\frac{k^{3}-k_{1}^{3}+k_{2}^{3}+k_{3}^{3}}{6}  }}{  k^{3}+(k_{1}+k_{3})^{3}+k_{2}^{3}}f_{k_{1}}g_{k_{2}}g_{k_{3}}.
\label{eq:fk-second}
\end{align}
and
\begin{align}
-\frac{i}{2}k\sum \limits_{k_{1}+k_{2}=k}^{\ast}e^{-it\Big(-k+k_{1}-k_{2}+\frac{k^{3}-k_{1}^{3}+k_{2}^{3}}{6}\Big)}&f_{k_{1}}g_{k_{2}} = \partial_{t}\Bigg(\sum \limits_{k_{1}+k_{2}=k}^{\ast} 3k \frac{e^{-it\frac{k^{3}-k_{1}^{3}+k_{2}^{3}}{6}   }}{ k^{3}-k_{1}^{3}+k_{2}^{3} }f_{k_{1}}g_{k_{2}}\Bigg)\notag\\
& +\frac{9i}{4} \sum \limits_{\substack{   k_{1}+k_{2}+k_{3}=k  \\   k_{1}+k_{2}  \neq 0 }}^{\ast} k(k_{1}+k_{2})\frac{e^{-it\frac{k^{3}-k_{1}^{3}-k_{2}^{3}+k_{3}^{3}}{6}  }}{  k^{3}-(k_{1}+k_{2})^{3}+k_{3}^{3} }f_{k_{1}}f_{k_{2}}g_{k_{3}}\notag\\
& -\frac{3i}{4} \sum \limits_{\substack{  k_{1}+k_{2}+k_{3}=k   \\   k_{2}+k_{3}  \neq 0 }}^{\ast} k(k_{2}+k_{3})\frac{e^{-it\frac{k^{3}+k_{1}^{3}+k_{2}^{3}+k_{3}^{3}}{6}   }}{ k^{3}+k_{1}^{3}-(k_{2}+k_{3})^{3} }g_{k_{1}}g_{k_{2}}g_{k_{3}}\notag\\
& +\frac{3i}{2} \sum \limits_{\substack{  k_{1}+k_{2}+k_{3}=k   \\   k_{1}+k_{3}  \neq 0 }}^{\ast} k(k_{1}+k_{3})\frac{e^{-it\frac{k^{3}-k_{1}^{3}+k_{2}^{3}+k_{3}^{3}}{6}   }}{  k^{3}-(k_{1}+k_{3})^{3}+k_{2}^{3}  }f_{k_{1}}g_{k_{2}}g_{k_{3}}\notag\\
& -\frac{3i}{4} \sum \limits_{\substack{  k_{1}+k_{2}+k_{3}=k   \\   k_{2}+k_{3}  \neq 0 }}^{\ast} k(k_{2}+k_{3})\frac{e^{-it\frac{k^{3}-k_{1}^{3}-k_{2}^{3}-k_{3}^{3}}{6}   }}{  k^{3}-k_{1}^{3}+(k_{2}+k_{3})^{3} }f_{k_{1}}f_{k_{2}}f_{k_{3}}\notag\\
& +\frac{9i}{4} \sum \limits_{\substack{  k_{1}+k_{2}+k_{3}=k   \\   k_{2}+k_{3}  \neq 0 }}^{\ast} k(k_{2}+k_{3})\frac{e ^{-it\frac{k^{3}-k_{1}^{3}+k_{2}^{3}+k_{3}^{3}}{6}  }}{ k^{3}-k_{1}^{3}+(k_{2}+k_{3})^{3}   }f_{k_{1}}g_{k_{2}}g_{k_{3}}\notag\\
& +\frac{3i}{2} \sum \limits_{\substack{   k_{1}+k_{2}+k_{3}=k  \\   k_{2}+k_{3}  \neq 0 }}^{\ast} k(k_{2}+k_{3})\frac{e ^{-it\frac{k^{3}-k_{1}^{3}-k_{2}^{3}+k_{3}^{3}}{6}  }}{  k^{3}-k_{1}^{3}+(k_{2}+k_{3})^{3}  }f_{k_{1}}f_{k_{2}}g_{k_{3}}.
\label{eq:fk-third}
\end{align}
  
We now analyze the resonant regions characterized by the relations
\begin{equation*}
(k_{1}+k_{2})(k_{1}+k_{3})( k_{2}+k_{3})  =0 \text{ and }(k_{1}+k_{3})( k_{2}+k_{3})  =0.
\end{equation*} 
We first evaluate the resonant contributions corresponding to terms $R_{1}$ and $R_{10}$.
Utilizing the conservation relation $k_{1}+k_{2}+k_{3}=k$, and the algebraic identities $k^{3}-k_{1}^{3}-(k_{2}+k_{3})^{3}=3kk_{1}(k_{2}+k_{3})$ and $k^{3}-k_{1}^{3}-k_{2}^{3}-k_{3}^{3}=3(k_{1}+k_{2})(k_{1}+k_{3})( k_{2}+k_{3})$,  we split the initial summation as follows
\begin{equation*}
\frac{27i}{4} \sum \limits_{\substack{  k_{1}+k_{2}+k_{3}=k   \\   k_{2}+k_{3}  \neq 0 }}^{\ast} k(k_{2}+k_{3})\frac{e^{-it\frac{k^{3}-k_{1}^{3}-k_{2}^{3}-k_{3}^{3}}{6}  }}{   k^{3}-k_{1}^{3}-(k_{2}+k_{3})^{3} }f_{k_{1}}f_{k_{2}}f_{k_{3}}=R_{1}(f,f,f)_{k}+R_{1}^{B}(f,f,f)_{k},
\end{equation*}
where the non-resonant constituent $R_{1}(f,f,f)_{k}$ is given by:
\begin{equation*}
R_{1}(f,f,f)_{k}=\frac{27i}{4} \sum \limits_{\substack{  k_{1}+k_{2}+k_{3}=k   \\  (k_{1}+k_{2})(k_{1}+k_{3})( k_{2}+k_{3})  \neq 0 }}^{\ast} k(k_{2}+k_{3})\frac{e^{-it\frac{k^{3}-k_{1}^{3}-k_{2}^{3}-k_{3}^{3}}{6}  }}{  k^{3}-k_{1}^{3}-(k_{2}+k_{3})^{3}  }f_{k_{1}}f_{k_{2}}f_{k_{3}},
\end{equation*}
and the boundary resonant contribution $R_{1}^{B}(f,f,f)_{k}$ is defined as:
\begin{equation*}
R_{1}^{B}(f,f,f)_{k}= \frac{27i}{4} \sum \limits_{\substack{  k_{1}+k_{2}+k_{3}=k   \\  (k_{1}+k_{2})(k_{1}+k_{3})  = 0 \\  k_{2}+k_{3}\neq 0}}^{\ast} k(k_{2}+k_{3})\frac{e^{-it\frac{k^{3}-k_{1}^{3}-k_{2}^{3}-k_{3}^{3}}{6}  }}{  k^{3}-k_{1}^{3}-(k_{2}+k_{3})^{3}  }f_{k_{1}}f_{k_{2}}f_{k_{3}}.
\end{equation*}
To evaluate the term $ R_{1}^{B}$, we partition the summation domain into three pairwise disjoint regions:

$S_{1}=\{ k_{1}+k_{2}=0, k_{1}+k_{3}=0  \text{ and } k_{2}+k_{3}\neq 0 \}=\{ k_{1}=-k, \quad k_{2}
=k,\quad k_{3}=k  \},$

$S_{2}=\{ k_{1}+k_{2}=0, k_{1}+k_{3}\neq 0   \text{ and } k_{2}+k_{3}\neq 0  \}=\{ k_{1}=j, \quad k_{2}=-j, \quad  k_{3}=k \text{ and } \vert j\vert \neq  \vert k\vert \},$

$S_{3}=\{ k_{1}+k_{2}\neq 0, k_{1}+k_{3}=0  \text{ and } k_{2}+k_{3}\neq 0   \}=\{k_{1}=j,\quad  k_{2}=k, \quad k_{3}=-j \text{ and } \vert j\vert \neq  \vert k\vert \}.$
Direct evaluation over these regions yields:
\begin{align}
R_{1}^{B}(f,f,f)_{k}&=\frac{9i}{4} \sum \limits_{\ell=1}^{3}\sum \limits_{S_{\ell} }^{\ast}\frac{e^{-\frac{1}{2}it(k_{1}+k_{2})(k_{1}+k_{3})( k_{2}+k_{3})}}{ k_{1} }f_{k_{1}}f_{k_{2}}f_{k_{3}}\notag\\
&=-\frac{9i}{4}\frac{\vert f_{k}\vert^{2}f_{k}}{k}+ \frac{9i}{2}f_{k}\sum \limits_{ \vert j\vert \neq  \vert k\vert }^{\ast}\frac{\vert f_{j}\vert^{2}}{j}=-\frac{9i}{4}\frac{\vert f_{k}\vert^{2}f_{k}}{k}.
\label{eq:RB-1}
\end{align}
In an analogous manner, we process the term associated with $R_{10}$ by setting
\begin{equation*}
-\frac{3i}{4}\sum \limits_{\substack{  k_{1}+k_{2}+k_{3}=k   \\   k_{2}+k_{3}  \neq 0 }}^{\ast} k(k_{2}+k_{3})\frac{e^{-it\frac{k^{3}-k_{1}^{3}-k_{2}^{3}-k_{3}^{3}}{6}   }}{   k^{3}-k_{1}^{3}+(k_{2}+k_{3})^{3}}f_{k_{1}}f_{k_{2}}f_{k_{3}}=R_{10}(f,f,f)_{k}+R_{10}^{B}(f,f,f)_{k},
\end{equation*}
where the nonresonant and boundary(resonant) contributions are given, respectively, by:
\begin{equation*}
R_{10}(f,f,f)_{k}=-\frac{3i}{4}\sum \limits_{\substack{  k_{1}+k_{2}+k_{3}=k   \\   (k_{1}+k_{2})(k_{1}+k_{3})(k_{2}+k_{3})  \neq 0 }}^{\ast} k(k_{2}+k_{3})\frac{e^{-it\frac{k^{3}-k_{1}^{3}-k_{2}^{3}-k_{3}^{3}}{6}   }}{    k^{3}-k_{1}^{3}+(k_{2}+k_{3})^{3}}f_{k_{1}}f_{k_{2}}f_{k_{3}},
\end{equation*}
and
\begin{equation*}
R_{10}^{B}(f,f,f)_{k}=-\frac{3i}{4}\sum \limits_{\substack{  k_{1}+k_{2}+k_{3}=k   \\ (k_{1}+k_{2})(k_{1}+k_{3})=0 \\ k_{2}+k_{3}  \neq 0 }}^{\ast} k(k_{2}+k_{3})\frac{e^{-it\frac{k^{3}-k_{1}^{3}-k_{2}^{3}-k_{3}^{3}}{6}   }}{   k^{3}-k_{1}^{3}+(k_{2}+k_{3})^{3}}f_{k_{1}}f_{k_{2}}f_{k_{3}}.
\end{equation*}
By applying the same disjoint partitions $S_{1}$, $S_{2}$ and $S_{3}$ as established for $R_{1}$, we find:
\begin{align}
R_{10}^{B}(f,f,f)_{k}& =-\frac{3i}{4}\sum \limits_{\ell=1}^{3}\sum \limits_{S_{\ell}}\frac{ke^{-\frac{1}{2}it(k_{1}+k_{2})(k_{1}+k_{3})( k_{2}+k_{3})}}{2k^{2}-kk_{1}+2k_{1}^{2}}f_{k_{1}}f_{k_{2}}f_{k_{3}}\notag\\
&=-\frac{3i}{20}\frac{\vert f_{k}\vert^{2}f_{k}}{k}-\frac{3i}{2}kf_{k}\sum \limits_{\vert j \vert \neq \vert k \vert}\frac{\vert f_{j}\vert^{2}}{2k^{2}-kj+2j^{2}}.
\label{eq:RB-10}
\end{align}
Combining (\ref{eq:RB-1}) and  (\ref{eq:RB-10}), the total boundary resonant contribution, denoted by $\mathcal{R}_{(f,f,f)_{k}}$,  is given by:
\begin{equation}
\label{eq:Rfff}
\mathcal{R}_{(f,f,f)_{k}}=-\frac{12i}{5}\frac{\vert f_{k}\vert^{2}f_{k}}{k}-\frac{3i}{2}kf_{k}\sum \limits_{\vert j \vert \neq \vert k \vert}\frac{\vert f_{j}\vert^{2}}{2k^{2}-kj+2j^{2}}.
\end{equation}
Following an argument analogous to the one used in the computation of the resonant term $\mathcal{R}_{(f,f,f)_{k}}$, we address
the resonant region \begin{equation*}
(k_{1}+k_{2})(k_{1}+k_{3})( k_{2}+k_{3})  =0,
\end{equation*} 
as it manifests in both $\mathcal{R}_{(f,f,g)_{k}}$ and $ \mathcal{R}_{(g,g,g)_{k}}$. To rigorously control the interaction $\mathcal{R}_{(f,f,g)_{k}}$, the resonant contribution paired with the $R_{3}$-term must be analyzed simultaneously with the corresponding resonant terms arising alongside $R_{4}$, $R_{7}$ and $R_{12}$. Similarly, isolating the resonant contribution associated to $R_{5}$ requires a joint examination with the resonant terms associated to $R_{8}$ in $ \mathcal{R}_{(g,g,g)_{k}}$. Consequently, we obtain 
\begin{align}
\mathcal{R}_{(f,f,g)_{k}} &={} -\frac{6i}{5}\frac{e^{-it\frac{k^{3}}{3}   }}{k}\vert f_{k}\vert^{2} g_{k} -\frac{3i}{8}\frac{e^{it\frac{k^{3}}{3}}}{k}f_{k}f_{k}g_{-k} +\frac{3i}{2}f_{k}\sum \limits_{\vert j \vert \neq \vert k \vert}^{\ast}\frac{e^{it\frac{j^{3}}{3}}}{j}f_{j}g_{-j}\notag\\
&\quad{}+\frac{3i}{2}f_{k}\sum \limits_{\vert j \vert \neq \vert k \vert}^{\ast}\frac{(j+k)e^{it\frac{j^{3}}{3}}}{2k^{2}+3jk+3j^{2}}f_{j}g_{-j}-\frac{9i}{2}kf_{k}\sum \limits_{\vert j \vert \neq \vert k \vert}^{\ast}\frac{(j+k)e^{it\frac{j^{3}}{3}   }}{j(2j^{2}+3kj+3k^{2})}f_{j}g_{-j}\notag\\
&\quad{}+ \frac{3i}{2}ke^{-it\frac{k^{3}}{3}}g_{k}\sum \limits_{\vert j \vert \neq \vert k \vert}^{\ast}\frac{\vert f_{j}\vert^{2}}{2k^{2}-kj+2j^{2}} +\frac{3i}{2}kf_{k}\sum \limits_{\vert j \vert \neq \vert k \vert}^{\ast}\frac{e^{it\frac{j^{3}}{3}}}{2k^{2}-kj+2j^{2}}f_{j}g_{-j},
\label{eq:Rffg}
\end{align}
\begin{align}
\mathcal{R}_{(g,g,g)_{k}} &={} -\frac{3i}{8}\frac{e^{-it\frac{k^{3}}{3}}}{k}\vert g_{k}\vert^{2}g_{k} -\frac{9i}{2}e^{-it\frac{k^{3}}{3}}g_{k}\sum \limits_{\vert j \vert \neq  \vert k \vert}^{\ast}\frac{k-j}{2k^{2}-3kj+3j^{2}}\vert g_{j}\vert^{2}\notag\\
&\quad{}+6ike^{-it\frac{k^{3}}{3}}g_{k}\sum \limits_{\substack{ j\geq 1\\ j\neq \vert k\vert}}^{\ast}\frac{j^{2}}{(2j^{2}-3kj+3k^{2})(2j^{2}+3kj+3k^{2})}\vert g_{j}\vert^{2}.
\label{eq:Rggg}
\end{align}
The resonant term $\mathcal{R}_{(f,g,g)_{k}}$, associated with the $R_{6}$ and $R_{9}$ terms, can be expressed by
\begin{align}
\mathcal{R}_{(f,g,g)_{k}} &={} -\frac{3i}{2}\sum \limits_{\vert j \vert \neq  \vert k \vert}^{\ast}\frac{k-j}{2k^{2}-3kj+3j^{2}}f_{k}\vert g_{j}\vert^{2}-\sum \limits_{\vert j \vert \neq  \vert k \vert}^{\ast}\frac{6ikf_{k}j^{2}}{(2j^{2}-3kj+3k^{2})(2j^{2}+3kj+3k^{2})}\vert g_{j}\vert^{2}\notag\\
&\quad{} -\frac{3i}{4}\frac{f_{k}\vert g_{k}\vert^{2}}{k}.
\label{eq:Rfgg}
\end{align}
Therefore, by combining (\ref{eq:fk-first}), (\ref{eq:fk-second}), (\ref{eq:fk-third}),  together with the resonant terms  (\ref{eq:Rfff}), \eqref{eq:Rffg}, \eqref{eq:Rggg} and \eqref{eq:Rfgg} we obtain the expression for $\partial_{t}f_{k}$ as
\begin{align*}
\partial_{t}f_{k}&=\partial_{t}\Bigg( \sum \limits_{k_{1}+k_{2}=k}^{\ast} \frac{9}{2}k \frac{e^{-it\frac{k^{3}-k_{1}^{3}-k_{2}^{3}}{6}   }}{  k^{3}-k_{1}^{3}-k_{2}^{3} }f_{k_{1}}f_{k_{2}}\Bigg)\\
        &+\frac{27i}{4} \sum \limits_{\substack{  k_{1}+k_{2}+k_{3}=k   \\   (k_{1}+k_{2})(k_{1}+k_{3})( k_{2}+k_{3})  \neq 0 }}^{\ast} k(k_{2}+k_{3})\frac{e^{-it\frac{k^{3}-k_{1}^{3}-k_{2}^{3}-k_{3}^{3}}{6}  }}{  k^{3}-k_{1}^{3}-(k_{2}+k_{3})^{3} }f_{k_{1}}f_{k_{2}}f_{k_{3}}\\
         &- \frac{9i}{4}\sum \limits_{\substack{  k_{1}+k_{2}+k_{3}=k   \\   k_{2}+k_{3}  \neq 0 }}^{\ast} k(k_{2}+k_{3})\frac{e^{-it\frac{k^{3}-k_{1}^{3}+k_{2}^{3}+k_{3}^{3}}{6}  }}{  k^{3}-k_{1}^{3}-(k_{2}+k_{3})^{3} }f_{k_{1}}g_{k_{2}}g_{k_{3}}\\
         &+\frac{9i}{2}  \sum \limits_{\substack{   k_{1}+k_{2}+k_{3}=k   \\   (k_{1}+k_{2})(k_{1}+k_{3})( k_{2}+k_{3})  \neq 0 }}^{\ast} k(k_{2} + k_{3})\frac{e^{-it\frac{k^{3}-k_{1}^{3}-k_{2}^{3}+k_{3}^{3}}{6}  }}{  k^{3}-k_{1}^{3}-(k_{2}+k_{3})^{3} }f_{k_{1}}f_{k_{2}}g_{k_{3}}\\
         &-\partial_{t}\Bigg(\sum \limits_{k_{1}+k_{2}=k}^{\ast} \frac{3k}{2} \frac{e^{-it\frac{k^{3}+k_{1}^{3}+k_{2}^{3}}{6}   }}{ k^{3}+k_{1}^{3}+k_{2}^{3}  }g_{k_{1}}g_{k_{2}}\Bigg)\\
        &+\frac{3i}{4} \sum \limits_{\substack{  k_{1}+k_{2}+k_{3}=k   \\    (k_{1}+k_{2})(k_{1}+k_{3})( k_{2}+k_{3}) \neq 0 }}^{\ast} k(k_{1}+k_{2})\frac{e^{-it\frac{k^{3}-k_{1}^{3}-k_{2}^{3}+k_{3}^{3}}{6}  }}{ k^{3}+(k_{1}+k_{2})^{3}+k_{3}^{3} }f_{k_{1}}f_{k_{2}}g_{k_{3}}\\
    &- \frac{9i}{4}\sum \limits_{\substack{  k_{1}+k_{2}+k_{3}=k   \\   (k_{1}+k_{2})(k_{1}+k_{3})( k_{2}+k_{3})  \neq 0 }}^{\ast} k(k_{2}+k_{3})\frac{e^{-it\frac{k^{3}+k_{1}^{3}+k_{2}^{3}+k_{3}^{3}}{6}   }}{  k^{3}+k_{1}^{3}+(k_{2}+k_{3})^{3} }g_{k_{1}}g_{k_{2}}g_{k_{3}}\\
    &-\frac{3i}{2}  \sum \limits_{\substack{  k_{1}+k_{2}+k_{3}=k   \\  (k_{1}+k_{3})( k_{2}+k_{3})  \neq 0 } }^{\ast}k(k_{1} + k_{3})\frac{e^{-it\frac{k^{3}-k_{1}^{3}+k_{2}^{3}+k_{3}^{3}}{6}  }}{  k^{3}+(k_{1}+k_{3})^{3}+k_{2}^{3}}f_{k_{1}}g_{k_{2}}g_{k_{3}}\\
    &+\partial_{t}\Bigg(\sum \limits_{k_{1}+k_{2}=k}^{\ast} 3k \frac{e^{-it\frac{k^{3}-k_{1}^{3}+k_{2}^{3}}{6}   }}{ k^{3}-k_{1}^{3}+k_{2}^{3} }f_{k_{1}}g_{k_{2}}\Bigg)\\
    &+\frac{9i}{4} \sum \limits_{\substack{   k_{1}+k_{2}+k_{3}=k  \\   (k_{1}+k_{2})(k_{1}+k_{3})( k_{2}+k_{3})  \neq 0 }}^{\ast} k(k_{1}+k_{2})\frac{e^{-it\frac{k^{3}-k_{1}^{3}-k_{2}^{3}+k_{3}^{3}}{6}  }}{  k^{3}-(k_{1}+k_{2})^{3}+k_{3}^{3} }f_{k_{1}}f_{k_{2}}g_{k_{3}}\\
    &-\frac{3i}{4} \sum \limits_{\substack{  k_{1}+k_{2}+k_{3}=k   \\  (k_{1}+k_{2})(k_{1}+k_{3})( k_{2}+k_{3})   \neq 0 }}^{\ast} k(k_{2}+k_{3})\frac{e^{-it\frac{k^{3}+k_{1}^{3}+k_{2}^{3}+k_{3}^{3}}{6}   }}{ k^{3}+k_{1}^{3}-(k_{2}+k_{3})^{3} }g_{k_{1}}g_{k_{2}}g_{k_{3}}\\
    &+\frac{3i}{2} \sum \limits_{\substack{  k_{1}+k_{2}+k_{3}=k   \\   (k_{1}+k_{3})( k_{2}+k_{3})  \neq 0 }}^{\ast} k(k_{1}+k_{3})\frac{e^{-it\frac{k^{3}-k_{1}^{3}+k_{2}^{3}+k_{3}^{3}}{6}   }}{  k^{3}-(k_{1}+k_{3})^{3}+k_{2}^{3}  }f_{k_{1}}g_{k_{2}}g_{k_{3}}\\
        &-\frac{3i}{4} \sum \limits_{\substack{  k_{1}+k_{2}+k_{3}=k   \\   (k_{1}+k_{2})(k_{1}+k_{3})( k_{2}+k_{3})  \neq 0 }}^{\ast} k(k_{2}+k_{3})\frac{e^{-it\frac{k^{3}-k_{1}^{3}-k_{2}^{3}-k_{3}^{3}}{6}   }}{  k^{3}-k_{1}^{3}+(k_{2}+k_{3})^{3} }f_{k_{1}}f_{k_{2}}f_{k_{3}}\\
        &+\frac{9i}{4} \sum \limits_{\substack{  k_{1}+k_{2}+k_{3}=k   \\   k_{2}+k_{3}  \neq 0 }}^{\ast} k(k_{2}+k_{3})\frac{e ^{-it\frac{k^{3}-k_{1}^{3}+k_{2}^{3}+k_{3}^{3}}{6}  }}{ k^{3}-k_{1}^{3}+(k_{2}+k_{3})^{3}   }f_{k_{1}}g_{k_{2}}g_{k_{3}}\\
        &+\frac{3i}{2} \sum \limits_{\substack{   k_{1}+k_{2}+k_{3}=k  \\  (k_{1}+k_{2})(k_{1}+k_{3})( k_{2}+k_{3})  \neq 0 }}^{\ast} k(k_{2}+k_{3})\frac{e ^{-it\frac{k^{3}-k_{1}^{3}-k_{2}^{3}+k_{3}^{3}}{6}  }}{  k^{3}-k_{1}^{3}+(k_{2}+k_{3})^{3}  }f_{k_{1}}f_{k_{2}}g_{k_{3}} \\
        &-\frac{12i}{5}\frac{\vert f_{k}\vert^{2}f_{k}}{k}-\frac{3i}{2}kf_{k}\sum \limits_{\vert j \vert \neq \vert k \vert}^{\ast}\frac{\vert f_{j}\vert^{2}}{2k^{2}-kj+2j^{2}}-\frac{6i}{5}\frac{e^{-it\frac{k^{3}}{3}   }}{k}\vert f_{k}\vert^{2} g_{k}\\
        & -\frac{3i}{8}\frac{e^{it\frac{k^{3}}{3}}}{k}f_{k}f_{k}g_{-k}+\frac{3i}{2}f_{k}\sum \limits_{\vert j \vert \neq \vert k \vert}^{\ast}\frac{e^{it\frac{j^{3}}{3}}}{j}f_{j}g_{-j} +\frac{3i}{2}f_{k}\sum \limits_{\vert j \vert \neq \vert k \vert}^{\ast}\frac{(j+k)e^{it\frac{j^{3}}{3}}}{2k^{2}+3jk+3j^{2}}f_{j}g_{-j}\\
       &-\frac{9i}{2}kf_{k}\sum \limits_{\vert j \vert \neq \vert k \vert}^{\ast}\frac{(j+k)e^{it\frac{j^{3}}{3}   }}{j(2j^{2}+3kj+3k^{2})}f_{j}g_{-j}+ \frac{3i}{2}ke^{-it\frac{k^{3}}{3}}g_{k}\sum \limits_{\vert j \vert \neq \vert k \vert}^{\ast}\frac{\vert f_{j}\vert^{2}}{2k^{2}-kj+2j^{2}}\\
       &+\frac{3i}{2}kf_{k}\sum \limits_{\vert j \vert \neq \vert k \vert}^{\ast}\frac{e^{it\frac{j^{3}}{3}}}{2k^{2}-kj+2j^{2}}f_{j}g_{-j}-\frac{3i}{8}\frac{e^{-it\frac{k^{3}}{3}}}{k}\vert g_{k}\vert^{2}g_{k}\\
       &-\frac{9i}{2}e^{-it\frac{k^{3}}{3}}g_{k}\sum \limits_{\vert j \vert \neq  \vert k \vert}^{\ast}\frac{k-j}{2k^{2}-3kj+3j^{2}}\vert g_{j}\vert^{2}+\sum \limits_{\substack{ j\geq 1\\ j\neq \vert k\vert}}^{\ast}\frac{6ike^{-it\frac{k^{3}}{3}}g_{k}j^{2}}{(2j^{2}-3kj+3k^{2})(2j^{2}+3kj+3k^{2})}\vert g_{j}\vert^{2} \\
      &-\frac{3i}{2}\sum \limits_{\vert j \vert \neq  \vert k \vert}^{\ast}\frac{f_{k}(k-j)}{2k^{2}-3kj+3j^{2}}\vert g_{j}\vert^{2}-\frac{3i}{4}\frac{f_{k}\vert g_{k}\vert^{2}}{k}-\hspace{-0.2cm}\sum \limits_{\vert j \vert \neq  \vert k \vert}^{\ast}\frac{6ikf_{k}j^{2}}{(2j^{2}-3kj+3k^{2})(2j^{2}+3kj+3k^{2})}\vert g_{j}\vert^{2}.
  \end{align*}
To obtain the desired result, we substitute  $f_{k}(t)=e^{-it\frac{k^{3}}{6}}u_{k}(t)$ and $g_{k}(t)=e^{it\frac{k^{3}}{6}}v_{k}(t)$ into the expression above.  
\end{proof}
\begin{proposition}\label{eq:v_k}
    Let $(u_{0},v_{0}) \in H^{s}_{0}(\mathbb{T})\times  H^{s}_{0}(\mathbb{T})$. Then we have 
\begin{align}
\partial_t \Bigg[
e^{it\frac{k^{3}}{6}}
\Big(
v_k
+\widetilde{B}_1(u,u)_k
-\widetilde{B}_2(v,v)_k
- &\widetilde{B}_3(u,v)_k
\Big)
\Bigg] =
e^{it\frac{k^{3}}{6}}
\Bigg[
\widetilde{R}_1(u,u,u)_k
+ \widetilde{R}_2(u,v,v)_k\notag\\
&+ \widetilde{R}_3(u,u,v)_k+ \widetilde{R}_4(u,u,v)_k
+ \widetilde{R}_5(v,v,v)_k\notag\\
&+ \widetilde{R}_6(u,v,v)_k+ \widetilde{R}_{7}(u,u,v)_k+ \widetilde{R}_{8}(v,v,v)_k\notag\\
&+ \widetilde{R}_{9}(u,v,v)_k
+ \widetilde{R}_{10}(u,u,u)_k+ \widetilde{R}_{11}(u,v,v)_k\notag\\
&+ \widetilde{R}_{12}(u,u,v)_k+ \widetilde{\mathcal{R}}_{(u,u,u)_{k}}+ \widetilde{\mathcal{R}}_{(u,u,v)_{k}}\notag\\
&+ \widetilde{\mathcal{R}}_{(u,v,v)_{k}}+ \widetilde{\mathcal{R}}_{(v,v,v)_{k}}\Bigg],
\label{eq:normal-second-eq}
\end{align}
where
\begin{gather*}
\widetilde{B}_1(u,u)_k=\sum \limits_{k_{1}+k_{2}=k}^{\ast}\frac{3k}{2}\frac{u_{k_{1}}u_{k_{2}}}{-k^{3}-k_{1}^{3}-k_{2}^{3}},\\
\widetilde{B}_2(v,v)_k=\sum \limits_{k_{1}+k_{2}=k}^{\ast}\frac{9k}{2}\frac{v_{k_{1}}v_{k_{2}}}{-k^{3}+k_{1}^{3}+k_{2}^{3}},\\
\widetilde{B}_3(u,v)_k= \sum \limits_{k_{1}+k_{2}=k}^{\ast}3k\frac{u_{k_{1}}v_{k_{2}}}{-k^{3}-k_{1}^{3}+k_{2}^{3}},\\
\widetilde{R}_1(u,u,u)_k=-\frac{9i}{4}\sum \limits_{\substack{  k_{1}+k_{2}+k_{3}=k   \\   (k_{1}+k_{2})(k_{1}+k_{3})(k_{2}+k_{3})  \neq 0 }}^{\ast}\frac{k(k_{2}+k_{3})u_{k_{1}}u_{k_{2}}u_{k_{3}}}{-k^{3}-k_{1}^{3}-(k_{2}+k_{3})^{3}},\\
\widetilde{R}_2(u,v,v)_k=\frac{3i}{4}\sum \limits_{\substack{   k_{1}+k_{2}+k_{3}=k  \\   (k_{1}+k_{2})(k_{1}+k_{3})(k_{2}+k_{3})  \neq 0 }}^{\ast}\frac{k(k_{2}+k_{3}) \, u_{k_{1}}v_{k_{2}}v_{k_{3}}}{-k^{3}-k_{1}^{3}-(k_{2}+k_{3})^{3}},\\
\widetilde{R}_3(u,u,v)_k=-\frac{3i}{2}\sum \limits_{\substack{  k_{1}+k_{2}+k_{3}=k   \\   (k_{1}+k_{2})(k_{2}+k_{3})  \neq 0 }}^{\ast}\frac{k(k_{2}+k_{3}) \, u_{k_{1}}u_{k_{2}}v_{k_{3}}}{-k^{3}-k_{1}^{3}-(k_{2}+k_{3})^{3}},\\
\widetilde{R}_4(u,u,v)_k=-\frac{9i}{4}\sum \limits_{\substack{  k_{1}+k_{2}+k_{3}=k   \\   k_{1}+k_{2}  \neq 0 }}^{\ast}\frac{k(k_{1}+k_{2})\, u_{k_{1}}u_{k_{2}}v_{k_{3}}}{-k^{3}+(k_{1}+k_{2})^{3}+k_{3}^{3}},\\
\widetilde{R}_5(v,v,v)_k=\frac{27i}{4}\sum \limits_{\substack{ k_{1}+k_{2}+k_{3}=k    \\   (k_{1}+k_{2})(k_{1}+k_{3})(k_{2}+k_{3})  \neq 0 }}^{\ast}\frac{k(k_{2}+k_{3})\, v_{k_{1}}v_{k_{2}}v_{k_{3}}}{-k^{3}+k_{1}^{3}+(k_{2}+k_{3})^{3}},\\
\widetilde{R}_6(u,v,v)_k=\frac{9i}{2}\sum \limits_{\substack{  k_{1}+k_{2}+k_{3}=k   \\  (k_{1}+k_{2})(k_{1}+k_{3})(k_{2}+k_{3})   \neq 0 }}^{\ast}\frac{k(k_{1}+k_{3})\, u_{k_{1}}v_{k_{2}}v_{k_{3}}}{-k^{3}+(k_{1}+k_{3})^{3}+k_{2}^{3}},\\
\widetilde{R}_{7}(u,u,v)_k=\frac{9i}{4} \! \sum \limits_{\substack{   k_{1}+k_{2}+k_{3}=k  \\   k_{1}+k_{2}  \neq 0 }}^{\ast}\frac{k(k_{1}+k_{2})\, u_{k_{1}}u_{k_{2}}v_{k_{3}}}{-k^{3}-(k_{1}+k_{2})^{3}+k_{3}^{3}},\\
\widetilde{R}_{8}(v,v,v)_k=-\frac{3i}{4} \! \sum \limits_{\substack{  k_{1}+k_{2}+k_{3}=k   \\   (k_{1}+k_{2})(k_{1}+k_{3})(k_{2}+k_{3})  \neq 0 }}^{\ast}\frac{k(k_{2}+k_{3})\, v_{k_{1}}v_{k_{2}}v_{k_{3}}}{-k^{3}+k_{1}^{3}-(k_{2}+k_{3})^{3}},\\
\widetilde{R}_{9}(u,v,v)_k=\frac{3i}{2} \! \! \sum \limits_{\substack{   k_{1}+k_{2}+k_{3}=k  \\  (k_{1}+k_{2})(k_{1}+k_{3})(k_{2}+k_{3})  \neq 0 }}^{\ast}\frac{k(k_{1}+k_{3})\, u_{k_{1}}v_{k_{2}}v_{k_{3}}}{-k^{3}-(k_{1}+k_{3})^{3}+k_{2}^{3}},\\
\widetilde{R}_{10}(u,u,u)_k=-\frac{3i}{4} \! \sum \limits_{\substack{  k_{1}+k_{2}+k_{3}=k   \\   (k_{1}+k_{2})(k_{1}+k_{3})(k_{2}+k_{3})  \neq 0 }}^{\ast}\frac{k(k_{2}+k_{3})\, u_{k_{1}}u_{k_{2}}u_{k_{3}}}{-k^{3}-k_{1}^{3}+(k_{2}+k_{3})^{3}},\\
\widetilde{R}_{11}(u,v,v)_k=\frac{9i}{4} \! \sum \limits_{\substack{  k_{1}+k_{2}+k_{3}=k   \\  (k_{1}+k_{2})(k_{1}+k_{3})(k_{2}+k_{3})  \neq 0 }}^{\ast}\frac{k(k_{2}+k_{3})\, u_{k_{1}}v_{k_{2}}v_{k_{3}}}{-k^{3}-k_{1}^{3}+(k_{2}+k_{3})^{3}},\\
\widetilde{R}_{12}(u,u,v)_k=\frac{3i}{2} \! \sum \limits_{\substack{  k_{1}+k_{2}+k_{3}=k   \\   (k_{1}+k_{2})(k_{2}+k_{3})  \neq 0 }}^{\ast}\frac{k(k_{2}+k_{3})\, u_{k_{1}}u_{k_{2}}v_{k_{3}}}{-k^{3}-k_{1}^{3}+(k_{2}+k_{3})^{3}},
\end{gather*}
\begin{align*}
\widetilde{\mathcal{R}}_{(u,u,u)_{k}} =&-6iku_{k}\sum \limits_{\substack{ j\geq 1\\ j\neq \vert k\vert}}^{\ast}\frac{j^{2}\vert u_{j}\vert^{2}}{(2j^{2}-3kj+3k^{2})(2j^{2}+3kj+3k^{2})}  +\frac{9i}{2}u_{k}\sum \limits_{\vert j \vert \neq  \vert k \vert}\frac{(k-j)\vert u_{j}\vert^{2}}{2k^{2}-3kj+3j^{2}}\\
&+\frac{3i}{8}\frac{\vert u_{k}\vert^{2}u_{k}}{k}       ,
\end{align*}
\begin{align*}
\widetilde{\mathcal{R}}_{(u,v,v)_{k}} &={} \frac{3i}{8}\frac{u_{-k}v_{k}v_{k}}{k} -\frac{3i}{2}v_{k}\sum \limits_{\vert j \vert \neq \vert k\vert}^{\ast}\frac{(k-j)u_{j}\overline{v_{j}}}{2k^{2}-3jk+3j^{2}} +\frac{6i}{5} \frac{u_{k}\vert v_{k}\vert^{2}}{k} +\frac{3i}{2}v_{k}\sum \limits_{\vert j \vert \neq \vert k\vert}^{\ast}\frac{u_{j}\overline{v_{j}}}{j}\\
&-\frac{3i}{2}v_{k}\sum \limits_{\vert j \vert \neq \vert k\vert}^{\ast}\frac{ ku_{j}\overline{v_{j}}}{2k^{2}+kj+2j^{2}} -\frac{3i}{2}u_{k}\sum \limits_{\vert j \vert \neq \vert k\vert}^{\ast}\frac{k \vert v_{j}\vert^{2}}{2k^{2}-kj+2j^{2}}\\
&\quad{}-\frac{9i}{2}v_{k}\sum \limits_{\vert j \vert \neq \vert k\vert}^{\ast}\frac{k(k-j)u_{j}\overline{v_{j}}}{j(2j^{2}-3kj+3k^{2})},
\end{align*}
\begin{align*}
\widetilde{\mathcal{R}}(u,u,v) &={} \frac{3i}{2}v_{k}\sum \limits_{\vert j \vert \neq \vert k \vert }^{\ast} \frac{(k-j)\vert u_{j}\vert^{2}}{2k^{2}-3kj+3j^{2}} +6iv_{k}\sum \limits_{\substack{j\neq  \vert k \vert\\ j\geq 1}}^{\ast} \frac{kj^{2}\vert u_{j}\vert^{2}}{(2j^{2}-3kj+3k^{2})(2j^{2}+3kj+3k^{2})}\\
&\quad{}+\frac{3i}{4}\frac{\vert u_{k}\vert^{2}v_{k}}{k},
\end{align*}
    \begin{align*}
\hspace{-3.5cm}\widetilde{\mathcal{R}}(v,v,v)_{k} &={} \frac{12i}{5}\frac{\vert v_{k}\vert^{2}v_{k}}{k} +\frac{3i}{2}v_{k}\sum \limits_{\vert j \vert \neq \vert k \vert}^{\ast}\frac{k\vert v_{j}\vert^{2}}{2k^{2}-kj+2j^{2}}.
\end{align*}
\end{proposition}
\begin{proof}
The proof follows by the same argument as in the proof of  Proposition \ref{eq:u_k}.
\end{proof}

\subsection{Estimation of bilinear and resonant terms}
\label{subsec:Bilinear-boundary}
We take the integral of the equations in  Proposition~\ref{eq:u_k} and Proposition~\ref{eq:v_k} from $0$ to $t$
\begin{align}
u_{k}-e^{it\frac{k^{3}}{6}}u_{k}(0)
&\quad=B_{1}(u,u)_{k}-B_{2}(v,v)_{k}+B_{3}(u,v)_{k}\notag\\
&\qquad+e^{it\frac{k^{3}}{6}}\Bigg[-B_{1}(u,u)_{k}(0)+B_{2}(v,v)_{k}(0)-B_{3}(u,v)_{k}(0)\Bigg]\notag\\
&\qquad+\int \limits_{0}^{t}e^{-i(t'-t)\frac{k^{3}}{6}}\Big(\sum \limits_{i=1}^{12}R_{i}+\mathcal{R}_{(u,u,u)}+\mathcal{R}_{(u,u,v)}+\mathcal{R}_{(u,v,v)}+ \mathcal{R}_{(v,v,v)}\Big)(t')dt' ,
\end{align}
\begin{align}
v_{k}-e^{-it\frac{k^{3}}{6}}v_{k}(0)
&\quad=-\widetilde{B}_{1}(u,u)_{k}+\widetilde{B}_{2}(v,v)_{k}+\widetilde{B}_{3}(u,v)_{k}\notag\\
&\qquad+e^{-it\frac{k^{3}}{6}}\Bigg[\widetilde{B}_{1}(u,u)_{k}(0)-\widetilde{B}_{2}(v,v)_{k}(0)-\widetilde{B}_{3}(u,v)_{k}(0)\Bigg]\notag\\
&\qquad+\int \limits_{0}^{t}e^{i(t'-t)\frac{k^{3}}{6}}\Big(\sum \limits_{i=1}^{12}\widetilde{R}_{i}+\widetilde{\mathcal{R}}_{(u,u,u)}+\widetilde{\mathcal{R}}_{(u,u,v)}+\widetilde{\mathcal{R}}_{(u,v,v)}+ \widetilde{\mathcal{R}}_{(v,v,v)}\Big)(t')dt'.
\end{align}
\begin{proposition}\label{prop:prop-B-1-2-3}
     For $s>-\frac{1}{2}$ and $a\leq 1$, we have
    \begin{equation*}
\Vert B_{1}(u,u)+ B_{2}(v,v)+ B_{3}(u,v)\Vert_{H^{s+a}}\lesssim \Big( \Vert u \Vert_{H^{s}}+\Vert v \Vert_{H^{s}}\Big)^{2}.
\end{equation*}
\end{proposition}
\begin{proof} Since the estimation of $ B_{1}(u,u)$ and $B_{2}(v,v)$ terms can be established by using  arguments analogous to those employed in the proof of $B_{3}(u,v)$, we present only the proof for the term $B_{3}(u,v)$.
    \begin{equation*}
\Vert B_{3}(u,v)\Vert_{H^{s+a}}\sim \Bigg\Vert \sum \limits_{k_{1}+k_{2}=k}^{\ast} \frac{\langle k \rangle^{s+a+1}u_{k_{1}}v_{k_{2}}}{k^{3}-k_{1}^{3}+k_{2}^{3}}\Bigg\Vert_{\ell_{k}^{2}}.
\end{equation*} 
    Observe that, as the denominator does not vanish, we have
    \begin{equation*}
\vert k^{3}-k_{1}^{3}+k_{2}^{3}  \vert =  \vert k_{2}(2k_{2}^{2}-3kk_{2}+3k^{2})  \vert \gtrsim \vert k_{2}\vert \max\{\vert k \vert , \vert k_{1} \vert, \vert k_{2} \vert \}^{2}.
\end{equation*}
    \textbf{Case A}: $\vert k_{1}\vert \geq \vert k_{2}\vert$, which implies that $\vert k_{1}\vert \gtrsim \vert k\vert$.
    \begin{align*}
\Vert B_{3}(u,v)\Vert_{H^{s+a}}
\lesssim \Bigg\Vert \sum \limits_{k_{1}+k_{2}=k}^{\ast} \frac{\langle k \rangle^{s+a+1}}{\langle k_{1}\rangle^{s+2}\langle k_{2}\rangle^{s+1}}\langle k_{1}\rangle^{s}u_{k_{1}}\langle k_{2}\rangle^{s}v_{k_{2}}\Bigg\Vert_{\ell_{k}^{2}}
&\lesssim  \Bigg\Vert  \langle k \rangle^{a-1}   \Big( \langle k \rangle^{s}u_{k} \ast \frac{\langle k\rangle^{s}v_{k}}{\langle k\rangle^{s+1}}\Big) \Bigg\Vert_{\ell_{k}^{2}}\\
&\lesssim \Vert u\Vert_{H^{s}}\Vert v\Vert_{H^{s}}
\end{align*}
    whenever $a\leq 1$ and $s> -\frac{1}{2}$. \\
\textbf{Case B}: $\vert k_{1}\vert < \vert k_{2}\vert$, which implies that $\vert k_{2}\vert \gtrsim \vert k\vert$.
    \begin{align*}
\Vert B_{3}(u,v)\Vert_{H^{s+a}}
\lesssim\Bigg\Vert\sum \limits_{k_{1}+k_{2}=k}^{\ast} \frac{ \langle k \rangle^{s+a+1}}{\langle k_{1}\rangle^{s+1}\langle k_{2}\rangle^{s+2}}\langle k_{1}\rangle^{s}u_{k_{1}}\langle k_{2}\rangle^{s}v_{k_{2}}\Bigg\Vert_{\ell_{k}^{2}}
&\lesssim \Bigg\Vert  \langle k \rangle^{a-1}   \Big(   \frac{\langle k\rangle^{s}u_{k}}{\langle k\rangle^{s+1}} \ast  \langle k \rangle^{s}v_{k} \Big) \Bigg\Vert_{\ell_{k}^{2}}\\
& \lesssim \Vert u\Vert_{H^{s}}\Vert v\Vert_{H^{s}}
\end{align*}
    if  $a\leq 1$ and $s>-\frac{1}{2}$. 
\end{proof}

\begin{proposition}
    For $s>-\frac{1}{2}$ and $a\leq 1$, we have
     \begin{equation*}
\Big\Vert  \widetilde{B}_{1}(u,u)+\widetilde{B}_{2}(v,v)+\widetilde{B}_{3}(u,v)\Big\Vert_{H^{s+a}}\lesssim \Big( \Vert u \Vert_{H^{s}}+\Vert v \Vert_{H^{s}}\Big)^{2}.
\end{equation*}
\end{proposition}
\begin{proof}
    The proof follows the same arguments  as those used in the proof of Proposition \ref{prop:prop-B-1-2-3}.
\end{proof}

\begin{proposition}\label{prop:Ruuv}
    For $s>- \frac{1}{2}$ and $a< \min\{1,2s+1\}$, we have
    \begin{equation*}
\Vert \mathcal{R}_{(u,u,v)} \Vert_{H^{s+a}} \lesssim \Vert u \Vert_{H^{s}}^{2}\Vert v \Vert_{H^{s}}.
\end{equation*}
\end{proposition}

\begin{proof}
First note that the first and second contributions to $\mathcal{R}_{(u,u,v)}$ are estimated in \cite{Erdogan-Tzirakis-global-sol-per} for $s>-\frac{1}{2}$ and $a<\min\{1,2s+1\}$. Therefore, we focus on estimating the third, fourth and fifth contributions to $\mathcal{R}_{(u,u,v)}$ as the remaining terms can be treated similarly.  We begin with the third and fifth terms. Observe that

\begin{equation*}
\Bigg \vert \frac{1}{j}-\frac{3k(j+k)}{j(2j^{2}+3kj+3k^{2})}\Bigg\vert  \lesssim  \frac{1}{\max\{\vert j\vert,\vert k \vert\}}.
\end{equation*}
    Hence, for $s\geq 0$
    \begin{align*}
\sum \limits_{\vert j \vert \neq \vert k \vert}^{\ast}\Bigg \vert \Big( \frac{1}{j}-\frac{3k(j+k)}{j(2j^{2}+3kj+3k^{2})}\Big)u_{j}\overline{v_{j}}\Bigg\vert
\lesssim  \sum \limits_{\vert j \vert \neq \vert k \vert}^{\ast}\frac{1}{\max\{\vert j\vert,\vert k \vert\}}\vert u_{j}\overline{v_{j}}\vert &\leq \frac{1}{\vert k \vert}\Vert u_{j} \Vert_{\ell^{2}_{j}}\Vert v_{j} \Vert_{\ell^{2}_{j}}\\
&\quad\leq \frac{1}{\vert k \vert}\Vert u \Vert_{H^{s}}\Vert v \Vert_{H^{s}}.
\end{align*} 
Consequently,
   \begin{align*}
\Bigg\Vert \mathcal{F}^{-1}\Bigg(\frac{3i}{2}u_{k}\sum \limits_{\vert j \vert \neq \vert k \vert}^{\ast} \Big( \frac{1}{j}-\frac{3k(j+k)}{j(2j^{2}+3kj+3k^{2})}\Big)u_{j}\overline{v_{j}}\Bigg)\Bigg\Vert_{H^{s+a}}
&\lesssim \Big\Vert \frac{\langle k\rangle^{s+a}}{\langle k\rangle}u_{k}\Big\Vert_{\ell_{k}^{2}}\Vert u\Vert_{H^{s}}\Vert v\Vert_{H^{s}}\\
&\quad\lesssim \Vert u\Vert_{H^{s}}^{2}\Vert v\Vert_{H^{s}}
\end{align*}
   whenever $a\leq 1$.\\
Now, let $-\frac{1}{2}<s<0$. Then
\begin{align*}
\sum \limits_{\vert j \vert \neq \vert k \vert}^{\ast}\left \vert \Big( \frac{1}{j}-\frac{3k(j+k)}{j(2j^{2}+3kj+3k^{2})}\Big)u_{j}\overline{v_{j}}\right\vert
&\leq \sum \limits_{\vert j \vert \neq \vert k \vert}^{\ast}\frac{1}{\max\{\vert j\vert^{2s+1},\vert k \vert^{2s+1}\}}\vert j \vert^{s}\vert u_{j}\vert \vert j \vert^{s} \vert v_{j}\vert\\
&\leq \frac{1}{\vert k \vert^{2s+1}}\Vert u \Vert_{H^{s}}\Vert v \Vert_{H^{s}}.
\end{align*} 
   Therefore,
 \begin{align*}
\Bigg\Vert \mathcal{F}^{-1}\Bigg(\frac{3i}{2}u_{k}\sum \limits_{\vert j \vert \neq \vert k \vert}^{\ast} \Big( \frac{1}{j}-\frac{3k(j+k)}{j(2j^{2}+3kj+3k^{2})}\Big)u_{j}\overline{v_{j}}\Bigg)\Bigg\Vert_{H^{s+a}}
&\lesssim \Big\Vert \frac{\langle k\rangle^{s+a}}{\langle k\rangle^{2s+1}}u_{k}\Big\Vert_{\ell_{k}^{2}}\Vert u\Vert_{H^{s}}\Vert v\Vert_{H^{s}}\\
&\quad\lesssim \Vert u\Vert_{H^{s}}^{2}\Vert v\Vert_{H^{s}}
\end{align*}
whenever $a\leq 2s+1$.

We next consider the fourth contribution to $\mathcal{R}_{(u,u,v)}$. Observe that   
  \begin{equation*}
\Bigg\vert \frac{j+k}{2k^{2}+3jk+3j^{2}}\Bigg\vert \lesssim  \frac{1}{\max\{\vert j\vert,\vert k\vert\}}.
\end{equation*}
Hence, arguing as above, we obtain
  \begin{align*}
\Bigg\Vert \mathcal{F}^{-1} \Bigg( \frac{3i}{2}u_{k}\sum \limits_{\vert j\vert \neq \vert k\vert }^{\ast}\frac{j+k}{2k^{2}+3jk+3j^{2}}u_{j}\overline{v_{j}}\Bigg)\Bigg\Vert_{H^{s+a}}
\lesssim \Vert u \Vert_{H^{s}}^{2}\Vert v \Vert_{H^{s}}
\end{align*}
  whenever $s> - \frac{1}{2}$ and $a\leq \min\{1,2s+1\}$.
\end{proof}

\begin{proposition}\label{prop:Rvvv}
     For $s>- \frac{1}{2}$ and $a< \min\{1,2s+1\}$, we have
    \begin{equation*}
\Vert \mathcal{R}_{(v,v,v)} \Vert_{H^{s+a}} \lesssim \Vert v \Vert_{H^{s}}^{3}.
\end{equation*}
\end{proposition}
\begin{proof}
   The first and second terms can be estimated analogously to those in Proposition  \ref{prop:Ruuv}. For the remaining term, we note that
    \begin{equation*}
\Bigg\vert\frac{j^{2}}{(2j^{2}-3kj+3k^{2})(2j^{2}+3kj+3k^{2})}\Bigg\vert\lesssim \frac{1}{\max\{ \vert j\vert^{2},\vert k \vert^{2}  \}}.
\end{equation*} 
   The desired conclusion then follows by repeating the argument used above.
\end{proof}
\begin{proposition}\label{prop:Ruuu-Ruvv}
     For $s>- \frac{1}{2}$ and $a< \min\{1,2s+1\}$, we have
    \begin{equation*}
\Vert \mathcal{R}_{(u,u,u)}+\mathcal{R}_{(u,v,v)} \Vert_{H^{s+a}} \lesssim \Vert u \Vert_{H^{s}} \Big(\Vert v \Vert_{H^{s}}^{2}+\Vert u \Vert_{H^{s}}^{2}\Big).
\end{equation*}
\end{proposition}

The proof is identical to the arguments used in Propositions \ref{prop:Rvvv}  and \ref{prop:Ruuv} and is therefore omitted.

\begin{proposition}\label{prop:mathcal-R-tilde}
    For $s>- \frac{1}{2}$ and $a< \min\{1,2s+1\}$, we have
     \begin{align*}
&\Big\Vert \widetilde{\mathcal{R}}_{(u,u,u)} + \widetilde{\mathcal{R}}_{(u,v,v)} +  \widetilde{\mathcal{R}}_{(u,u,v)} + \widetilde{\mathcal{R}}_{(v,v,v)} \Big\Vert_{H^{s+a}}\\*
&\quad\lesssim \big( \Vert u \Vert_{H^{s}}+  \Vert v \Vert_{H^{s}}\big)^{3}.
\end{align*}
\end{proposition}
\begin{proof}
    The proof proceeds by adapting the techniques developed in Propositions   \ref{prop:Ruuv}, \ref{prop:Rvvv} and \ref{prop:Ruuu-Ruvv}. 
\end{proof}

\subsection{Estimation of nonresonant terms}
\label{subsec:Nonresonant}
The nonlinear terms $R_{3}$, $R_{4}$, $R_{7}$ and $R_{12}$, corresponding to the $(u,u,v)$ interactions, exhibit degeneracies in their associated phase functions in certain frequency regions, leading to a loss of frequency decay in the denominators. However, for the terms $R_{4}$ and $R_{12}$, additional factors in the denominators compensate for this loss: no such compensating factors are available for the terms $R_{3}$ and $R_{7}$. Therefore,  we present only the estimates for the $R_{3}$ and  $R_{4}$ terms for this interaction, as the remaining estimates follow by a similar argument.

Let $\delta > 0$ denote the local existence time, and define $\psi_{\delta}(t) = \psi\!\left(\tfrac{t}{\delta}\right)$, where $\psi$ is a compactly supported Schwartz function on $[-2,2]$ such that $\psi \equiv 1$ on $[-1,1]$ and decays rapidly on $[-2,-1) \cup (1,2]$. By the embedding $X^{s+a,b} \hookrightarrow L^{\infty}_{t \in [-\delta,\delta]} H^{s+a}_x:=L_{t}^{\infty}H_{x}^{s+a}$ together with Lemma~\ref{eq:generic-lemma}, we deduce 
\begin{align*}
&\Bigg\Vert u-e^{-t\frac{\partial_{xxx}}{6}}u_{0}\Bigg\Vert_{L_{t}^{\infty}H_{x}^{s+a}}\\
&\lesssim \Big\Vert\sum\limits_{i=1}^{3}B_{i}\Big\Vert_{L_{t}^{\infty}H_{x}^{s+a}}
\hspace{-0.1cm}+\left\Vert\int\limits_{0}^{t} e^{(t'-t)\frac{\partial_{xxx}}{6}}
 \Big(\sum\limits_{i=1}^{12}R_{i}\hspace{-0.1cm}+\hspace{-0.1cm}\mathcal{R}_{(u,u,u)}\hspace{-0.1cm}+\hspace{-0.1cm}\mathcal{R}_{(u,u,v)}
\hspace{-0.1cm}+\hspace{-0.1cm}\mathcal{R}_{(v,v,v)}\hspace{-0.1cm}+\hspace{-0.1cm}\mathcal{R}_{(u,v,v)}\Big)(t')dt'
 \right\Vert_{L_{t}^{\infty}H_{x}^{s+a}}\\
&\lesssim \Big\Vert\sum\limits_{i=1}^{3}B_{i}\Big\Vert_{L_{t}^{\infty}H_{x}^{s+a}}
\hspace{-0.1cm}+\hspace{-0.1cm}\Bigg\Vert\psi_{\delta}(t)\int\limits_{0}^{t}
 e^{(t'-t)\frac{\partial_{xxx}}{6}}\sum\limits_{i=1}^{12}R_{i}(t')dt'
 \Bigg\Vert_{X_{1/6,\delta}^{s+a,b}}
\hspace{-0.1cm}\\
&\hspace{7cm}+\Big\Vert\mathcal{R}_{(u,u,u)}+\mathcal{R}_{(u,u,v)}
 +\mathcal{R}_{(v,v,v)}+\mathcal{R}_{(u,v,v)}\Big\Vert_{L_{t}^{\infty}H_{x}^{s+a}}\\
&\lesssim \Big\Vert\sum\limits_{i=1}^{3}B_{i}\Big\Vert_{L_{t}^{\infty}H_{x}^{s+a}}
 +\Big\Vert\sum\limits_{i=1}^{12}R_{i}\Big\Vert_{X_{1/6,\delta}^{s+a,b-1}}
+\Big\Vert\mathcal{R}_{(u,u,u)}+\mathcal{R}_{(u,u,v)}
 +\mathcal{R}_{(v,v,v)}+\mathcal{R}_{(u,v,v)}\Big\Vert_{L_{t}^{\infty}H_{x}^{s+a}}.
\end{align*}
and
\begin{align*}
\Bigg\Vert v-&e^{t\frac{\partial_{xxx}}{6}}v_{0}\Bigg\Vert_{L_{t}^{\infty}H_{x}^{s+a}}\\
&\lesssim \Big\Vert\sum\limits_{i=1}^{3}\widetilde{B}_{i}\Big\Vert_{L_{t}^{\infty}H_{x}^{s+a}}
 +\Big\Vert\sum\limits_{i=1}^{12}\widetilde{R}_{i}\Big\Vert_{X_{1/6,\delta}^{s+a,b-1}}
+\Big\Vert\widetilde{\mathcal{R}}_{(u,u,u)}+\widetilde{\mathcal{R}}_{(u,u,v)}
 +\widetilde{\mathcal{R}}_{(v,v,v)}+\widetilde{\mathcal{R}}_{(u,v,v)}
 \Big\Vert_{L_{t}^{\infty}H_{x}^{s+a}}.
\end{align*}
\begin{proposition}\label{eq:burak-tzirakis-R1}
    For $s>-\frac{1}{2}$ and $a< \min\{ 1, 2s+1\}$, we have
    \begin{equation*}
\Vert R_{1}(u,v,w)\Vert_{X_{1/6,\delta}^{s+a,b-1}}\lesssim \Vert u\Vert_{X_{1/6,\delta}^{s,1/2}}\Vert v\Vert_{X_{1/6,\delta}^{s,1/2}}\Vert w\Vert_{X_{1/6,\delta}^{s,1/2}}.
\end{equation*}
\end{proposition}
This proposition was established by Erdoğan and Tzirakis in \cite{Erdogan-Tzirakis-global-sol-per}.
\begin{proposition}\label{prop:R-2}
    For $s>-\frac{1}{2}$ and $a< \min\{ 1, 2s+1\}$, we have
    \begin{equation*}
\Vert R_{2}(u,v,w)\Vert_{X_{1/6,\delta}^{s+a,b-1}}\lesssim \Vert u\Vert_{X_{1/6,\delta}^{s,1/2}}\Vert v\Vert_{X_{-1/6,\delta}^{s,1/2}}\Vert w\Vert_{X_{-1/6,\delta}^{s,1/2}}.
\end{equation*}
\end{proposition}
\begin{proof}
By duality, it suffices to prove that
\begin{align*}
\Vert R_{2}(u,v,w)\Vert_{ X^{s+a,b-1}_{1/6,\delta}}&=\sup \limits_{g\in X_{1/6,\,\delta}^{\,-s-a,\;1-b}}\Bigg| \sum \limits_{k}\int \limits_{\R}\widehat{R_{2}}(k,\tau)\widehat{g}(-k,-\tau)d\tau \Bigg|\\
&\lesssim \Vert u \Vert_{ X^{s,1/2}_{1/6,\delta}}\Vert v \Vert_{ X^{s,1/2}_{-1/6,\delta}}\Vert w \Vert_{ X^{s,1/2}_{-1/6,\delta}} \Vert g \Vert_{ X^{-s-a,1-b}_{1/6,\delta}}.
\end{align*}
\begin{align*}
\widehat{R}_2(k,\tau)&=-\frac{9i}{4}\int \limits_{\tau_{1}+\tau_{2}+\tau_{3}=\tau}\sum \limits_{\substack{ k_{1}+k_{2}+k_{3}=k    \\   k_{2}+k_{3}  \neq 0 }}^{\ast}\frac{ k(k_{2}+k_{3})\, \widehat{u}(k_{1},\tau_{1})\widehat{v}(k_{2},\tau_{2})\widehat{w}(k_{3},\tau_{3})}{k^{3}-k_{1}^{3}-(k_{2}+k_{3})^{3}}\\
&=-\frac{3i}{4}\int \limits_{\tau_{1}+\tau_{2}+\tau_{3}=\tau}\sum \limits_{\substack{ k_{1}+k_{2}+k_{3}=k    \\   k_{2}+k_{3}  \neq 0 }}^{\ast}  \frac{\widehat{u}(k_{1},\tau_{1})\widehat{v}(k_{2},\tau_{2})\widehat{w}(k_{3},\tau_{3})}{k_{1}}.
\end{align*}
Define 
\begin{align*}
&f_{1}(k,\tau)=\big\langle k \big\rangle^{s}\big\langle \tau- \phi^{1/6}(k)\big\rangle^{\frac{1}{2}}\vert \widehat{u}(k,\tau)\vert,\\
&f_{2}(k,\tau)=\big\langle k \big\rangle^{s}\big\langle \tau- \phi^{-1/6}(k)\big\rangle^{\frac{1}{2}}\vert \widehat{v}(k,\tau)\vert,\\
&f_{3}(k,\tau)=\big\langle k \big\rangle^{s}\big\langle \tau- \phi^{-1/6}(k)\big\rangle^{\frac{1}{2}}\vert \widehat{w}(k,\tau)\vert,\\
&f_{4}(k,\tau)=\big\langle k \big\rangle^{-s-a}\big\langle \tau- \phi^{1/6}(k)\big\rangle^{1-b}\vert \widehat{g}(k,\tau)\vert.
\end{align*}
Then
\begin{align*}
\Vert R_{2}(u,v,w)\Vert_{ X^{s+a,b-1}_{1/6,\delta}}
\leq \hspace{-0.3cm}\sup \limits_{g\in X_{1/6,\,\delta}^{\,-s-a,\;1-b}} \hspace{-0.3cm}\mathop{\int \sum} \limits_{\substack{\tau_1+\tau_2+\tau_3+\tau_{4}=0 \\ k_1+k_2+k_3+k_{4}=0 \\k_{2}+k_{3}\neq 0}}\frac{\prod \limits_{i=1}^{4}f_{i}(k_{i},\tau_{i})}{\langle k_{1}\rangle^{s+1}\langle k_{2}\rangle^{s}\langle k_{3}\rangle^{s}\langle k_{4}\rangle^{-s-a}\prod \limits_{i=1}^{4}\big\langle \tau_{i}- \phi^{\alpha_{i}}(k_{i})\big \rangle^{1-b}},
\end{align*}
where $\alpha_{1}=\alpha_{4}=\frac{1}{6}$ and $\alpha_{2}=\alpha_{3}=-\frac{1}{6}$.\\
Observe the following  inequality
\begin{align*}
\langle \Phi\rangle&=\Big\langle \sum \limits_{i=1}^{4}\phi^{\alpha_{i}}(k_{i}) \Big\rangle \sim \Big|(k_{2}+k_{3})\big[2(k_{2}+k_{3})^{2}-3(k_{1}k_{4}+k_{2}k_{3})   \big]\Big|\\
&\sim \Big\vert (k_{2}+k_{3})\Big[ 2(k_{2}+k_{3})^{2}+3(k_{2}-k_{3})^{2}+3(2k_{1}+k_{2}+k_{3})^{2}\Big]\Big\vert\\
&\gtrsim \big|k_{2}+k_{3}\big|\vert k_{i}\vert^{2} \text{ for }i=1,2,3.
\end{align*}
Moreover, as $\tau_1+\tau_2+\tau_3+\tau_{4}=0$, it follows that 
\begin{align}
\max\limits_{i\in\{1,2,3,4\}}
 \big\langle\tau_{i}-\phi^{\alpha_{i}}(k_{i})\big\rangle
 :=\big\langle\tau_{i_{0}}-\phi^{\alpha_{i_{0}}}(k_{i_{0}})\big\rangle
 \gtrsim\langle\Phi\rangle
 \gtrsim\big|k_{2}+k_{3}\big|M^{2} \text{  for some }i_{0}\in\{1,2,3,4\},
\label{eq:max-tau-phi-ki}
\end{align}
where $
 M=\max\{\vert k_{1}\vert,\vert k_{2}\vert,\vert k_{3}\vert\}.$
\begin{align}
\prod \limits_{i=1}^{4}\big\langle \tau_{i}- \phi^{\alpha_{i}}(k_{i})\big \rangle^{1-b}&=\big\langle \tau_{i_{0}}-\phi^{\alpha_{i_{0}}}(k_{i_{0}})\big\rangle^{1-b}\prod \limits_{\substack{i=1 \\i \neq i_{0} }}^{4}\big\langle \tau_{i}- \phi^{\alpha_{i}}(k_{i})\big \rangle^{1-2b}\prod \limits_{\substack{i=1 \\i \neq i_{0} }}^{4}\big\langle \tau_{i}- \phi^{\alpha_{i}}(k_{i})\big \rangle^{b}\notag\\
&\gtrsim\big\langle \tau_{i_{0}}-\phi^{\alpha_{i_{0}}}(k_{i_{0}})\big\rangle^{4-7b}\prod \limits_{\substack{i=1 \\i \neq i_{0} }}^{4}\big\langle \tau_{i}- \phi^{\alpha_{i}}(k_{i})\big \rangle^{b}.
\label{eq:tau-phi-separate}
\end{align}
Let $b=\frac{1}{2}+\epsilon$ for small $\epsilon>0$. We claim that 
\begin{equation}
\label{eq:I}
I=\frac{\vert k_{1}k_{2}k_{3}k_{4}\vert^{\epsilon}}{\vert k_{1} \vert^{s+1}\vert k_{2} \vert^{s}\vert k_{3} \vert^{s}\vert k_{4} \vert^{-a-s}\vert k_{2}+k_{3} \vert^{4-7b}M^{8-14b}}\lesssim 1.
\end{equation}
Suppose that the claim holds. By using ( \ref{eq:tau-phi-separate}), (\ref{eq:max-tau-phi-ki}) and (\ref{eq:I}), consecutively, we obtain the following
\begin{align*}
\mathop{\int \sum} \limits_{\substack{\tau_1+\tau_2+\tau_3+\tau_{4}=0 \\ k_1+k_2+k_3+k_{4}=0 \\k_{2}+k_{3}\neq 0}}&\frac{\prod \limits_{i=1}^{4}f_{i}(k_{i},\tau_{i})}{\langle k_{1}\rangle^{s+1}\langle k_{2}\rangle^{s}\langle k_{3}\rangle^{s}\langle k_{4}\rangle^{-s-a}\prod \limits_{i=1}^{4}\big\langle \tau_{i}- \phi^{\alpha_{i}}(k_{i})\big \rangle^{1-b}}\\
&\lesssim  \mathop{\int \sum} \limits_{\substack{\tau_1+\tau_2+\tau_3+\tau_{4}=0 \\ k_1+k_2+k_3+k_{4}=0 \\k_{2}+k_{3}\neq 0}}\frac{\prod \limits_{i=1}^{4}f_{i}(k_{i},\tau_{i})}{\prod \limits_{\substack{i=1 \\i \neq i_{0} }}^{4}\langle k_{i}\rangle^{\epsilon}\big\langle \tau_{i}- \phi^{\alpha_{i}}(k_{i})\big \rangle^{b}}.
\end{align*}
Without loss of generality, suppose that $i_{0}=i_{1}$ and define
\begin{equation*}
F_{1}(k_{1},\tau_{1}):=f_{1}(k_{1},\tau_{1}) \text{ and } F_{i}(k_{i},\tau_{i}):=\frac{f_{i}(k_{i},\tau_{i})}{\vert k_{i}\vert^{\epsilon}\big\langle \tau_{i}-\phi^{\alpha_{i}}(k_{i})\big \rangle^{b}} \text{ for }i=2,3,4.
\end{equation*}
Then
\begin{align*}
\mathop{\int \sum} \limits_{\substack{\tau_1+\tau_2+\tau_3+\tau_{4}=0 \\ k_1+k_2+k_3+k_{4}=0 \\k_{2}+k_{3}\neq 0}}&\frac{\prod \limits_{i=1}^{4}f_{i}(k_{i},\tau_{i})}{\langle k_{1}\rangle^{s+1}\langle k_{2}\rangle^{s}\langle k_{3}\rangle^{s}\langle k_{4}\rangle^{-s-a}\prod \limits_{i=1}^{4}\big\langle \tau_{i}- \phi^{\alpha_{i}}(k_{i})\big \rangle^{1-b}}\\
&\lesssim \Big( F_{1} \ast  F_{2} \ast  F_{3} \ast F_{4} \Big)(0,0)\\
&=\int \limits_{\R \times \mathbb{T}} \check{f}_{1}(x,t)\prod \limits_{i=1}^{4}\Bigg(  \frac{f_{i}(k,\tau)}{\vert k\vert^{\epsilon}\big\langle \tau-\phi^{\alpha_{i}}(k)\big \rangle^{b}} \Bigg)^{\vee}\\
&\lesssim \Vert f_{1}\Vert_{L^{2}}\prod \limits_{i=2}^{4} \Bigg \Vert \Bigg(  \frac{f_{i}(k,\tau)}{\vert k\vert^{\epsilon}\big\langle \tau-\phi^{\alpha_{i}}(k)\big \rangle^{b}} \Bigg)^{\vee} \Bigg \Vert_{L^{6}}\\
&\lesssim \prod \limits_{i=1}^{4}\Vert f_{i}\vert_{L^{2}}.
\end{align*}
In the last step, we used Proposition~\ref{eq:inverse-fourier-L6-norm}. It remains to prove the claim \eqref{eq:I}.
For the region $s>-\frac{1}{3}$,
we may assume that $a+s\geq 0$.
In this region, we have 
\begin{equation*}
\vert k_{1}k_{2}k_{3}\vert^{\epsilon-s}\lesssim M^{-\min{(s-\epsilon,3s-3\epsilon)}} \text{ and } \quad \vert k_{4}\vert^{a+s+\epsilon}\lesssim M^{a+s+\epsilon}.
\end{equation*}
By symmetry, it suffices to consider the following two cases:\par\nopagebreak[4]

\noindent\textbf{Case 1:} $M=\vert k_{1}\vert$.\par\nopagebreak[4]

\begin{equation*}
I\lesssim \frac{M^{-\min{(s-\epsilon,3s-3\epsilon)}}M^{a+s+\epsilon}}{M^{2-14\epsilon}}\lesssim 1 \quad \text{ if } a<\min(2,2+2s).
\end{equation*}

\noindent\textbf{Case 2:} $M=\vert k_{2}\vert$.\par\nopagebreak[4]

\begin{equation*}
I\lesssim \frac{M^{-\min{(s-\epsilon,3s-3\epsilon)}}M^{a+s+\epsilon}}{M^{1-14\epsilon}}\lesssim 1 \quad \text{ if } a<\min(1,1+2s).
\end{equation*}
For the region $-\frac{1}{2}<s\leq -\frac{1}{3}$, \\
we have that $a+s<3s+1\leq 0$.
Using the inequality $\vert k_{1}\vert ^{\frac{1}{2}-7\epsilon}\vert k_{2}+k_{3}\vert ^{\frac{1}{2}-7\epsilon}\gtrsim \vert k_{4}\vert ^{\frac{1}{2}-7\epsilon}$, we obtain
\begin{equation*}
I\lesssim \frac{\vert k_{1}\vert^{\epsilon}}{\vert k_{1}\vert^{\frac{1}{2}+7\epsilon+s}}\frac{\vert k_{2}\vert^{\epsilon-s}}{ M^{\frac{1}{2}-7\epsilon}}\frac{\vert k_{3}\vert^{\epsilon-s}}{ M^{\frac{1}{2}-7\epsilon}}\frac{\vert k_{4}\vert^{\epsilon}}{\vert k_{4}\vert^{-a-s+\frac{1}{2}-7\epsilon}}\lesssim 1
\end{equation*}
by choosing $\epsilon$ small enough.
\end{proof}
\begin{proposition} \label{prop:R-3}
  For $s> -\frac{1}{2}$ and $a< \min\{ 1,\frac{4s}{3}+\frac{2}{3}\}$, we have
    \begin{equation*}
\Vert R_{3}(u,v,w)\Vert_{X_{1/6,\delta}^{s+a,b-1}}\lesssim \Vert u\Vert_{X_{1/6,\delta}^{s,1/2}}\Vert v\Vert_{X_{1/6,\delta}^{s,1/2}}\Vert w\Vert_{X_{-1/6,\delta}^{s,1/2}}.
\end{equation*}   
\end{proposition}
\begin{proof}
 Following the same argument as in Proposition \ref{prop:R-2}, it suffices to prove that  
\begin{equation}
\label{eq:I-R-3}
I=\frac{\vert k_{1}k_{2}k_{3}k_{4}\vert^{\epsilon}}{\vert k_{1} \vert^{s+1}\vert k_{2} \vert^{s}\vert k_{3} \vert^{s}\vert k_{4} \vert^{-a-s}\langle \Phi \rangle^{\frac{1}{2}-7\epsilon}}\lesssim 1,
\end{equation}
where $\Phi=\frac{-2k_{3}^{3}-3(k_{1}+k_{2})(k_{1}+k_{3})(k_{2}+k_{3})}{6}$ with $\alpha_{1}=\alpha_{2}=\alpha_{4}=\frac{1}{6}$ and $\alpha_{3}=-\frac{1}{6}$.

For the proof of \eqref{eq:I-R-3}, we employ a low-high decomposition. Let $M=\max \limits_{i\in \{1,2,3\}}\vert k_{i}\vert$.
We now proceed case by case.\par\nopagebreak[4]

\noindent\textbf{Case 1:} $\vert k_{1}\vert \sim  \vert k_{2}\vert \sim \vert k_{3}\vert \sim M$. Then\par\nopagebreak[4]

\begin{equation*}
I=\frac{\vert k_{1}k_{2}k_{3}k_{4} \vert^{\epsilon}}{\vert k_{1}\vert^{s+1}\vert k_{2}\vert^{s}\vert k_{3}\vert^{s}\vert k_{4}\vert^{-a-s}\langle \Phi \rangle^{\frac{1}{2}-7\epsilon}}\lesssim \frac{M^{4\epsilon}}{M^{3s+1}\vert k_{4}\vert^{-a-s}\langle \Phi \rangle^{\frac{1}{2}-7\epsilon}}.
\end{equation*}

\noindent\textbf{Case 1A:} $\vert k_{4}\vert \sim M$.\par\nopagebreak[4]

\begin{equation*}
I\lesssim \frac{M^{4\epsilon}}{M^{2s+1-a}}\lesssim 1, \quad \text{ provided that }a<2s+1.
\end{equation*}

\noindent\textbf{Case 1B:} $\vert k_{4}\vert \ll M$. We decompose\par\nopagebreak[4]
\begin{equation*}
k_{1}=aM+L_{1}, \,k_{2}=bM+L_{2} \text{ and } k_{3}=-(a+b)M+L_{3},
\end{equation*}
where $L_{1},L_{2}$ and $L_{3}$ represent lower-order terms relative to $M$, while $a,b$ and $a+b$ are nonzero constants of order one.
It follows that \begin{equation*}
\Phi=\frac{(a+b)(2a^{2}+ab+2b^{2})M^{3}+L}{6},
\end{equation*} where $L$ is a lower-order term.  Since the coefficient of $M^{3}$ is non-zero and of order one, we conclude that \begin{equation*}
\langle \Phi \rangle \sim M^{3}.
\end{equation*}
If $s\geq 0$, then
\begin{equation*}
I\lesssim \frac{M^{4\epsilon+a} \vert k_{4}\vert ^{s}}{M^{3s+\frac{5}{2}-21\epsilon}}\lesssim 1, \quad \text{ provided that }a <2s+\frac{5}{2}.
\end{equation*}
If $-\frac{1}{2}<s<0$, then
\begin{equation*}
I\lesssim \frac{M^{4\epsilon+a}}{ M^{3s+\frac{5}{2}-21\epsilon}\vert k_{4}\vert ^{-s}}\lesssim  1, \quad \text{ provided that }a <3s+\frac{5}{2}.
\end{equation*}

\noindent\textbf{Case 2:} $\vert k_{1}\vert \ll  \vert k_{2}\vert \sim \vert k_{3}\vert \sim M$. Then\par\nopagebreak[4]
\begin{equation*}
I\lesssim \frac{M^{4\epsilon}}{\vert k_{1}\vert^{s+1}M^{2s}\vert k_{4}\vert^{-a-s}\langle \Phi \rangle^{\frac{1}{2}-7\epsilon}}.
\end{equation*}
 Accordingly, we write
\begin{equation*}
k_{1}= L_{1}, \,k_{2}=aM+L_{2}  \text{ and } k_{3}=bM+L_{3},
\end{equation*}
where $L_{1},L_{2}$ and $L_{3}$ represent lower-order terms relative to $M$, and $a,b$ are nonzero constants of order one.
It follows that 
\begin{equation*}
\Phi=\frac{-b(2b^{2}+3ab+3a^{2})M^{3}+L}{6},
\end{equation*} where $L$ is a lower-order term.  Since the coefficient of $M^{3}$ is nonzero and similar to one, we have 
\begin{equation*}
\langle \Phi \rangle \sim M^{3} \text{ and }I\lesssim \frac{M^{4\epsilon+a}}{M^{2s+\frac{3}{2}-21\epsilon}\vert k_{1}\vert^{s+1}\vert k_{4}\vert^{-s}}.
\end{equation*}
We conclude that \\
If $s\geq 0$, then 
\begin{equation*}
I\lesssim \frac{M^{4\epsilon +a+s}}{M^{2s+\frac{3}{2}-21\epsilon}}\lesssim 1, \quad \text{ provided that }a<s+\frac{3}{2}.
\end{equation*}
If $-\frac{1}{2}<s< 0$, then 
\begin{equation*}
I\lesssim \frac{M^{4\epsilon +a}}{M^{2s+\frac{3}{2}-21\epsilon}}\lesssim 1, \quad \text{ provided that }a<2s+\frac{3}{2}.
\end{equation*}

\noindent\textbf{Case 3:} $\vert k_{2}\vert \ll  \vert k_{1}\vert \sim \vert k_{3}\vert \sim M$.\\
Applying the same frequency decompositions to  $k_{1},$ $k_{2}$ and $k_{3}$ as in Case 2, we obtain
\begin{equation*}
\langle \Phi \rangle \sim M^{3} \text{ and }I\lesssim \frac{M^{4\epsilon+a}}{M^{2s+\frac{5}{2}-21\epsilon}\vert k_{2}\vert^{s}\vert k_{4}\vert^{-s}}.
\end{equation*}
If $s\geq 0$, then
\begin{equation*}
I\lesssim \frac{M^{4\epsilon +a+s}}{M^{2s+\frac{5}{2}-21\epsilon}}\lesssim 1, \quad \text{ provided that }a<s+\frac{5}{2}.
\end{equation*}
If $-\frac{1}{2}<s<0$, then
\begin{equation*}
I\lesssim \frac{M^{4\epsilon +a-s}}{M^{2s+\frac{5}{2}-21\epsilon}}\lesssim 1, \quad \text{ provided that }a<3s+\frac{5}{2}.
\end{equation*}

\noindent\textbf{Case 4:} $\vert k_{3}\vert \ll  \vert k_{1}\vert \sim \vert k_{2}\vert \sim M$.\par\nopagebreak[4]

\noindent\textbf{Case 4A:} $\vert 2k_{3}^{3}\vert \ll \vert 3(k_{1}+k_{2})(k_{1}+k_{3})(k_{2}+k_{3}) \vert$. Then\par\nopagebreak[4]
$\langle \Phi \rangle \sim \vert k_{1}+k_{2}\vert M^{2}.$
Hence, 
\begin{equation*}
I\lesssim \frac{M^{4\epsilon+a}}{M^{2s+2-14\epsilon}\vert k_{3}\vert^{s}\vert k_{4}\vert^{-s}\vert k_{3}+k_{4}\vert^{\frac{1}{2}-7\epsilon}}.
\end{equation*}
 By applying the inequalities $\vert k_{3}\vert \vert k_{3}+k_{4}\vert \gtrsim \vert k_{4}\vert $ and $\vert k_{4}\vert \vert k_{3}+k_{4}\vert \gtrsim \vert k_{3}\vert $, we obtain \\
if $s\geq \frac{1}{2}$, then
\begin{equation*}
I\lesssim \frac{M^{11\epsilon+a}\vert k_{4}\vert^{s-\frac{1}{2}}}{M^{2s+2-14\epsilon}}\lesssim \frac{M^{11\epsilon+a+s-\frac{1}{2}}}{M^{2s+2-14\epsilon}} \lesssim 1, \quad \text{ provided that  }a<s+\frac{5}{2},
\end{equation*}
if $0\leq s<\frac{1}{2}$, then
\begin{equation*}
I\lesssim \frac{M^{11\epsilon+a}\vert k_{4}\vert^{s}}{M^{2s+2-14\epsilon}\vert k_{3}\vert^{s} \vert k_{3}+k_{4}\vert^{s}}\lesssim \frac{M^{11\epsilon+a}}{M^{2s+2-14\epsilon}} \lesssim 1, \quad \text{ provided that  }a<2s+2,
\end{equation*}
if $-\frac{1}{2}<s<0$, 
\begin{equation*}
I\lesssim \frac{M^{11\epsilon+a}\vert k_{3}\vert^{-s}}{M^{2s+2-14\epsilon}\vert k_{4}\vert^{-s} \vert k_{3}+k_{4}\vert^{-s}}\lesssim \frac{M^{11\epsilon+a}}{M^{2s+2-14\epsilon}} \lesssim 1, \quad \text{ provided that  }a<2s+2.
\end{equation*}

\noindent\textbf{Case 4B:} $\vert 2k_{3}^{3}\vert \gg \vert 3(k_{1}+k_{2})(k_{1}+k_{3})(k_{2}+k_{3}) \vert$. In this region, we have $\langle \Phi \rangle \gg  \vert k_{1}+k_{2}\vert M^{2} $ and\par\nopagebreak[4]
$\vert k_{3}\vert\gg \vert k_{1}+k_{2}\vert^{\frac{1}{3}} M^{\frac{2}{3}}$. Furthermore, $\vert k_{1}+k_{2}\vert \ll \vert k_{3}\vert $, and hence $\vert k_{4} \vert \sim \vert k_{3}\vert$.
\begin{equation*}
I\lesssim \frac{M^{4\epsilon}}{M^{2s+2-14\epsilon}\vert k_{3}\vert^{-a}\vert k_{1}+k_{2}\vert^{\frac{1}{2}-7\epsilon}}\lesssim 1,  \quad \text{ provided that  }a<2s+2.
\end{equation*}

\noindent\textbf{Case 4C:} $\vert 2k_{3}^{3}\vert \sim \vert 3(k_{1}+k_{2})(k_{1}+k_{3})(k_{2}+k_{3}) \vert$.\par\nopagebreak[4]
In this region, we have 
$\vert k_{3}\vert^{3}\sim \vert k_{1}+k_{2}\vert M^{2}$,  and consequently $\vert k_{1}+k_{2}\vert \ll \vert k_{3}\vert$. It follows that $\vert k_{4}\vert  \sim \vert k_{3}\vert$. 

\begin{equation*}
I\lesssim \frac{M^{4\epsilon}}{M^{2s+1}\vert k_{4}\vert^{-a}}\lesssim 1, \text{ provided that }a<2s+1.
\end{equation*}

\noindent\textbf{Case 5:} $\vert k_{1} \vert , \vert k_{3}\vert \ll  \vert k_{2}\vert \sim  M$. Then we have that $\vert k_{4}\vert\sim M$.\par\nopagebreak[4]

\noindent\textbf{Case 5A:} $\vert 2k_{3}^{3}\vert \ll \vert 3(k_{1}+k_{2})(k_{1}+k_{3})(k_{2}+k_{3}) \vert$. Then\par\nopagebreak[4]
\begin{equation*}
\langle \Phi\rangle \sim \vert k_{1}+k_{3}\vert M^{2}.
\end{equation*} By using the inequality $\vert k_{1}\vert \vert k_{1}+k_{3}\vert \gtrsim \vert k_{3}\vert $, it follows that
\begin{equation*}
I\lesssim  \frac{M^{11\epsilon }}{M^{1-a-14\epsilon}\vert k_{1}\vert ^{s+1}\vert k_{3}\vert ^{s}\vert k_{1}+k_{3}\vert^{\frac{1}{2}}}\lesssim \frac{M^{11\epsilon }}{M^{1-a-14\epsilon}\vert k_{1}\vert^{s+\frac{1}{2}}\vert k_{3}\vert^{s+\frac{1}{2}}} \lesssim 1, \quad \text{ provided that } a<1.
\end{equation*}

\noindent\textbf{Case 5B:} $\vert 2k_{3}^{3}\vert \gg \vert 3(k_{1}+k_{2})(k_{1}+k_{3})(k_{2}+k_{3}) \vert$. In this region, we have\par\nopagebreak[4]
$\langle \Phi \rangle \gg \vert k_{1}+k_{3}\vert M^{2} $ and   $\vert k_{3} \vert \gtrsim \vert k_{1}+k_{3}\vert^{\frac{1}{3}} M^{\frac{2}{3}}$. Furthermore, $\vert k_{1}+k_{3}\vert \ll \vert k_{3}\vert $, thus $\vert k_{1}\vert \sim \vert k_{3}\vert $.
 
\begin{equation*}
I\lesssim \frac{M^{4\epsilon}}{M^{1-a-14\epsilon}\vert k_{3}\vert^{2s+1} \vert k_{1}+k_{3}\vert^{\frac{1}{2}-7\epsilon}}\lesssim 1, \quad \text{ provided that }a<\frac{4s}{3}+\frac{5}{3}.
\end{equation*}

\noindent\textbf{Case 5C:} $\vert 2k_{3}^{3}\vert \sim \vert 3(k_{1}+k_{2})(k_{1}+k_{3})(k_{2}+k_{3}) \vert$. Then $\vert k_{3} \vert \sim \vert k_{1}+k_{3}\vert^{\frac{1}{3}} M^{\frac{2}{3}}$. Similarly, $\vert k_{1}+k_{3}\vert \ll \vert k_{3}\vert $, and hence $\vert k_{1}\vert \sim \vert k_{3}\vert $.\par\nopagebreak[4]

\begin{equation*}
I\lesssim \frac{M^{4\epsilon}}{M^{-a}\vert k_{3}\vert^{2s+1}}\sim \frac{M^{4\epsilon}}{M^{-a}\vert k_{1}+k_{3}\vert^{\frac{2s+1}{3}}M^{\frac{4s+2}{3}}}\lesssim 1, \quad \text{ provided that }a<\frac{4s}{3}+\frac{2}{3}.
\end{equation*}

\noindent\textbf{Case 6:}  $\vert k_{1} \vert , \vert k_{2}\vert \ll  \vert k_{3}\vert \sim  M$. Then we have\par\nopagebreak[4]
$\vert k_{4}\vert \sim M$ and $\langle \Phi \rangle \sim M^{3}$.
\begin{equation*}
I\lesssim \frac{M^{4\epsilon}}{\vert k_{1}\vert^{s+1}\vert k_{2}\vert^{s}M^{-a+\frac{3}{2}-21\epsilon}}.
\end{equation*}
If  $s\geq 0$, then
\begin{equation*}
I\lesssim \frac{M^{4\epsilon}}{M^{-a+\frac{3}{2}-21\epsilon}}\lesssim 1 \quad \text{ provided that }a<\frac{3}{2}.
\end{equation*}
If $-\frac{1}{2}<s<0$, then
\begin{equation*}
I\lesssim \frac{M^{4\epsilon-s}}{M^{-a+\frac{3}{2}-21\epsilon}}\lesssim 1 \quad \text{ provided that }a<s+\frac{3}{2}.
\end{equation*}

\noindent\textbf{Case 7:}  $\vert k_{2} \vert , \vert k_{3}\vert \ll  \vert k_{1}\vert \sim  M$. Then we have $\vert k_{4} \vert \sim M$.\par\nopagebreak[4]

\noindent\textbf{Case 7A:} $\vert 2k_{3}^{3}\vert \ll \vert 3(k_{1}+k_{2})(k_{1}+k_{3})(k_{2}+k_{3}) \vert$.\par\nopagebreak[4]
\begin{equation*}
\vert \Phi \vert \sim \vert k_{2}+k_{3}\vert M^{2} \text{; consequently, } I\lesssim \frac{M^{4\epsilon}}{M^{2-a-14\epsilon}\vert k_{2} \vert^{s}\vert k_{3}\vert^{s}\vert k_{2}+k_{3}\vert^{\frac{1}{2}-7\epsilon}}.
\end{equation*}
If $s\geq 0,$ then
\begin{equation*}
I\lesssim \frac{M^{4\epsilon}}{M^{2-a-14\epsilon}}\lesssim 1 \quad \text{ provided that } a<2.
\end{equation*}
If $-\frac{1}{2}<s< 0,$  then
\begin{equation*}
I\lesssim \frac{M^{4\epsilon-2s}}{M^{2-a-14\epsilon}}\lesssim 1 \quad \text{ provided that } a<2s+2.
\end{equation*}

\noindent\textbf{Case 7B:} $\vert 2k_{3}^{3}\vert \gg \vert 3(k_{1}+k_{2})(k_{1}+k_{3})(k_{2}+k_{3}) \vert$. In this region, we have  $\langle \Phi \rangle\gg \vert k_{2}+k_{3}\vert M^{2} $ and  $\vert k_{3}\vert \gg \vert k_{2}+k_{3}\vert^{\frac{1}{3}} M^{\frac{2}{3}}$. Moreover, $\vert k_{2}+k_{3}\vert \ll \vert k_{3}\vert $, and hence $\vert k_{2}\vert \sim \vert k_{3}\vert  $.\par\nopagebreak[4]
 
\begin{equation*}
I\lesssim \frac{M^{4\epsilon}}{M^{2-a-14\epsilon}\vert k_{3}\vert^{2s}\vert k_{2}+k_{3}\vert^{\frac{1}{2}-7\epsilon} }.
\end{equation*}
If $s\geq 0,$  then
\begin{equation*}
I\lesssim \frac{M^{11\epsilon}}{M^{2+\frac{4s}{3}-a-14 \epsilon}\vert k_{2}+k_{3}\vert^{\frac{2s}{3}+\frac{1}{2}} }\lesssim 1 \quad \text{ provided that } a<\frac{4s}{3}+2.
\end{equation*}
If $-\frac{1}{2}<s< 0,$ then
\begin{equation*}
I\lesssim \frac{M^{4\epsilon}\vert k_{3}\vert^{-2s}}{M^{2-a-14\epsilon}\vert k_{2}+k_{3}\vert^{\frac{1}{2}-7\epsilon}}\lesssim 1 \quad \text{ provided that } a<2s+2.
\end{equation*}

\noindent\textbf{Case 7C:} $\vert 2k_{3}^{3}\vert \sim \vert 3(k_{1}+k_{2})(k_{1}+k_{3})(k_{2}+k_{3}) \vert$. Then  $\vert k_{3}\vert \sim \vert k_{2}+k_{3}\vert^{\frac{1}{3}} M^{\frac{2}{3}}$. Similarly, $\vert k_{2}+k_{3}\vert \ll \vert k_{3}\vert $, so $\vert k_{2}\vert \sim \vert k_{3}\vert  $.\par\nopagebreak[4]
\begin{equation*}
I\lesssim \frac{M^{4\epsilon}}{M^{1+\frac{4s}{3}-a}\vert k_{2}+k_{3}\vert^{\frac{2s}{3}}}.
\end{equation*}
If $s\geq 0$, then
\begin{equation*}
I\lesssim \frac{M^{4\epsilon}}{M^{1+\frac{4s}{3}-a}}\lesssim 1, \quad \text{ provided that } a<\frac{4s}{3}+1.
\end{equation*}
If $-\frac{1}{2}<s< 0,$ then
\begin{equation*}
I\lesssim \frac{M^{4\epsilon-\frac{2s}{3}}}{M^{1+\frac{4s}{3}-a}}\lesssim 1, \quad \text{ provided that } a<2s+1.
\end{equation*}

\end{proof}

\begin{proposition}\label{prop:R-4}
  For $s>-\frac{1}{2} $ and $a< \min\{2,2s+1,\frac{4s}{3}+1\}$, we have
    \begin{equation*}
\Vert R_{4}(u,v,w)\Vert_{X_{1/6,\delta}^{s+a,b-1}}\lesssim \Vert u\Vert_{X_{1/6,\delta}^{s,1/2}}\Vert v\Vert_{X_{1/6,\delta}^{s,1/2}}\Vert w\Vert_{X_{-1/6,\delta}^{s,1/2}}.
\end{equation*}    
\end{proposition}
\begin{proof}
  By proceeding as in the proof of Proposition \ref{prop:R-2}, it suffices to establish that
\begin{equation*}
I= \frac{\vert k_{1}+k_{2}\vert\vert k_{1}k_{2}k_{3}k_{4}\vert^{\epsilon}}{\prod \limits_{i=1}^{3}\langle k_{i}\rangle^{s}\langle k_{4}\rangle^{-s-a}\max\{\vert k_{3}\vert^{2} ,\vert k_{4}\vert^{2} \}\langle \Phi \rangle^{\frac{1}{2}-7\epsilon}}\lesssim 1,
\end{equation*}
where $\Phi=\sum \limits_{i=1}^{4}\phi^{\alpha_{i}}(k_{i})=\frac{-2k_{3}^{3}-3(k_{1}+k_{2})(k_{1}+k_{3})(k_{2}+k_{3})}{6}$ with $\alpha_{1}=\alpha_{2}=\alpha_{4}=\frac{1}{6}$ and $\alpha_{3}=-\frac{1}{6}$.

Let $M=\max \limits_{i\in \{1,2,3\}}\vert k_{i}\vert$. By symmetry, we may assume that $\vert k_{1} \vert \leq \vert k_{2}\vert $. We proceed by employing a low-high frequency decomposition.

\noindent\textbf{Case 1:} $\vert k_{1} \vert \sim  \vert k_{2} \vert \sim \vert k_{3} \vert \sim M$.\par\nopagebreak[4]

\noindent\textbf{Case 1A:} $\vert k_{4}\vert \sim M.$\par\nopagebreak[4]
\begin{equation*}
I\lesssim \frac{\vert k_{1}+k_{2}\vert M^{4\epsilon}}{M^{2s-a+2}}\lesssim \frac{ M^{4\epsilon+1}}{M^{2s-a+2}}\lesssim 1, \quad \text{ provided that } a<2s+1.
\end{equation*}

\noindent\textbf{Case 1B:} $\vert k_{4} \vert \ll M.$ By using the decomposition as in the proof of  Proposition \ref{prop:R-3}, we observe that $\langle \Phi \rangle \sim M^{3}$. Then\par\nopagebreak[4]
\begin{equation*}
I\lesssim \frac{M^{4\epsilon+a}}{M^{3s+\frac{5}{2}-21\epsilon}\vert k_{4}\vert^{-s}}.
\end{equation*}
If $s\geq 0$, then
\begin{equation*}
I\lesssim \frac{M^{4\epsilon+a+s}}{M^{3s+\frac{5}{2}-21\epsilon}}\lesssim 1, \quad \text{ provided that }a<2s+\frac{5}{2}.
\end{equation*}
If $-\frac{1}{2}<s< 0$, then
\begin{equation*}
I\lesssim \frac{M^{4\epsilon+a}}{M^{3s+\frac{5}{2}-21\epsilon}}\lesssim 1, \quad \text{ provided that }a<3s+\frac{5}{2}.
\end{equation*}

\noindent\textbf{Case 2:} $\vert k_{1} \vert \ll  \vert k_{2} \vert \sim \vert k_{3} \vert \sim M$. By the same argument used in Proposition \ref{prop:R-3}, it follows that $\langle \Phi \rangle \sim M^{3}$. Then\par\nopagebreak[4]
\begin{equation*}
I\lesssim \frac{M^{4\epsilon+a}}{M^{2s+\frac{5}{2}-21\epsilon}  \vert k_{1}\vert^{s}\vert k_{4}\vert^{-s}}.
\end{equation*}
If $s\geq 0$, then
\begin{equation*}
I\lesssim \frac{M^{4\epsilon+a+s}}{M^{2s+\frac{5}{2}-21\epsilon}\vert k_{1}\vert^{s}}\lesssim 1, \quad \text{ provided that }a<s+\frac{5}{2}.
\end{equation*}
If $-\frac{1}{2}<s< 0$, then
\begin{equation*}
I\lesssim \frac{M^{4\epsilon+a-s}}{M^{2s+\frac{5}{2}-21\epsilon}\vert k_{4}\vert^{-s}}\lesssim 1, \quad \text{ provided that }a<3s+\frac{5}{2}.
\end{equation*}

\noindent\textbf{Case 3:} $\vert k_{3} \vert \ll  \vert k_{1} \vert \sim \vert k_{2} \vert \sim M$.\par\nopagebreak[4]

\begin{equation*}
I\lesssim \frac{\vert k_{1}+k_{2}\vert M^{4\epsilon+a}}{M^{2s}\vert k_{3}\vert^{s}\vert k_{4}\vert^{-s}\max\{ \vert k_{3}\vert^{2},\vert k_{4}\vert^{2}\}\langle \Phi \rangle^{\frac{1}{2}-7\epsilon}}.
\end{equation*}

\noindent\textbf{Case 3A:} $\vert 2k_{3}^{3}\vert \ll \vert 3(k_{1}+k_{2})(k_{1}+k_{3})(k_{2}+k_{3}) \vert$. Then $\langle \Phi \rangle \sim \vert k_{1}+k_{2}\vert M^{2}$.\par\nopagebreak[4]
Thus
\begin{equation*}
I\lesssim \frac{M^{11\epsilon+a}}{M^{2s+1-14\epsilon}\vert k_{3}\vert^{s}\vert k_{4}\vert^{-s}\max\{\vert k_{3}\vert^{\frac{3}{2}} ,\vert k_{4}\vert^{\frac{3}{2}}\}}.
\end{equation*}
If $s\geq \frac{3}{2}$, then
\begin{equation*}
I\lesssim \frac{M^{11\epsilon+a}\vert k_{4}\vert^{s-\frac{3}{2}}}{M^{2s+1-14\epsilon}\vert k_{3}\vert^{s}}\lesssim \frac{M^{11\epsilon+a+s-\frac{3}{2}}}{M^{2s+1-14\epsilon}}\lesssim 1, \quad \text{ provided that }a<s+\frac{5}{2}.
\end{equation*}
If $0\leq s< \frac{3}{2}$, then
\begin{equation*}
I\lesssim \frac{M^{11\epsilon+a}}{M^{2s+1-14\epsilon}\vert k_{3}\vert^{s}\vert k_{4}\vert^{\frac{3}{2}-s}}\lesssim 1, \quad \text{ provided that }a<2s+1.
\end{equation*}
If $-\frac{1}{2}< s< 0$, then
\begin{equation*}
I\lesssim \frac{M^{11\epsilon+a}}{M^{2s+1-14\epsilon}\vert k_{3}\vert^{s+\frac{3}{2}}\vert k_{4}\vert^{-s}}\lesssim 1, \quad \text{ provided that }a<2s+1.
\end{equation*}

\noindent\textbf{Case 3B:} $\vert 2k_{3}^{3}\vert \gg \vert 3(k_{1}+k_{2})(k_{1}+k_{3})(k_{2}+k_{3}) \vert$. In this region, we have that  $\langle \Phi \rangle \gg \vert k_{1}+k_{2}\vert M^{2}$ and $\vert k_{3}\vert \gg \vert k_{1}+k_{2}\vert^{\frac{1}{3}}M^{\frac{2}{3}}  $. Furthermore, we have $\vert k_{1}+k_{2}\vert \ll \vert k_{3}\vert $, so it follows that $\vert k_{4}\vert \sim \vert k_{3} \vert $. \par \nopagebreak[4]

\begin{equation*}
I\lesssim \frac{\vert k_{1}+k_{2}\vert^{\frac{1}{2}}M^{11\epsilon+a}}{M^{2s+1-14\epsilon}\vert k_{3}\vert^{2}}\lesssim \frac{M^{11\epsilon+a}}{M^{2s+\frac{7}{3}-14\epsilon}\vert k_{1}+k_{2}\vert^{\frac{1}{6}}}\lesssim 1,\quad \text{ provided that }a<2s+\frac{7}{3}.
\end{equation*}

\noindent\textbf{Case 3C:} $\vert 2k_{3}^{3}\vert \sim \vert 3(k_{1}+k_{2})(k_{1}+k_{3})(k_{2}+k_{3}) \vert$. In this region,  we have that $\vert k_{3}\vert^{3}\sim \vert k_{1}+k_{2}\vert M^{2}$. Similarly, we have $\vert k_{1}+k_{2}\vert \ll \vert k_{3}\vert $, thus  $\vert k_{4}\vert \sim \vert k_{3} \vert $.\par\nopagebreak[4]

\begin{equation*}
I\lesssim \frac{\vert k_{1}+k_{2}\vert M^{4\epsilon+a}}{M^{2s}\vert k_{3}\vert^{2}}\lesssim  \frac{\vert k_{1}+k_{2}\vert^{\frac{1}{3}} M^{4\epsilon+a}}{M^{2s+\frac{4}{3}}}\lesssim 1,\quad \text{ provided that }a<2s+1.
\end{equation*}

\noindent\textbf{Case 4:} $\vert k_{1} \vert,  \vert k_{3} \vert \ll \vert k_{2} \vert \sim M$. Then\par\nopagebreak[4]
\begin{equation*}
\vert k_{4}\vert \sim M  \text{ and } I\lesssim \frac{M^{4\epsilon}}{M^{1-a}\vert k_{1}\vert^{s}\vert k_{3} \vert^{s}\langle \Phi \rangle^{\frac{1}{2}-7\epsilon}}.
\end{equation*}

\noindent\textbf{Case 4A:} $\vert 2k_{3}^{3}\vert \ll \vert 3(k_{1}+k_{2})(k_{1}+k_{3})(k_{2}+k_{3}) \vert$. Then $\langle \Phi \rangle \sim \vert k_{1}+k_{3}\vert M^{2} $.\par\nopagebreak[4]
\begin{equation*}
I\lesssim \frac{M^{4\epsilon}}{M^{2-a-14\epsilon}\vert k_{1}\vert^{s}\vert k_{3}\vert^{s}\vert k_{1}+k_{3}\vert^{\frac{1}{2}-7\epsilon}}.
\end{equation*}
If $s\geq 0$, then
\begin{equation*}
I\lesssim \frac{M^{4\epsilon}}{M^{2-a-14\epsilon}}\lesssim 1,\quad \text{ provided that }a<2.
\end{equation*}
If $-\frac{1}{2}<s<0$, then
\begin{equation*}
I\lesssim \frac{M^{4\epsilon}\vert k_{1}\vert^{-s} \vert k_{3}\vert^{-s}}{M^{2-a-14\epsilon}\vert k_{1}+k_{3}\vert^{\frac{1}{2}-7\epsilon}}\lesssim \frac{M^{4\epsilon-2s}}{M^{2-a-14\epsilon}}\lesssim 1,\quad \text{ provided that }a<2s+2.
\end{equation*}

\noindent\textbf{Case 4B:} $\vert 2k_{3}^{3}\vert \gg \vert 3(k_{1}+k_{2})(k_{1}+k_{3})(k_{2}+k_{3}) \vert$. Then  we have that $\langle \Phi \rangle \gg \vert k_{1}+k_{3}\vert M^{2} $ and $\vert k_{3}\vert \gg \vert k_{1}+k_{3}\vert^{\frac{1}{3}}M^{\frac{2}{3}}$. Moreover, we have $\vert k_{1}+k_{3}\vert \ll \vert k_{3}\vert $, and consequently $\vert k_{1}\vert \sim \vert k_{3}\vert $.\par\nopagebreak[4]
\begin{equation*}
I\lesssim \frac{M^{4\epsilon}}{M^{2-a-14\epsilon}\vert k_{3}\vert^{2s}\vert k_{1}+k_{3}\vert^{\frac{1}{2}-7\epsilon}}.
\end{equation*}
If $s\geq 0$, then
\begin{equation*}
I\lesssim \frac{M^{4\epsilon}}{M^{2-a+\frac{4s}{3}-14\epsilon}\vert k_{1}+k_{3}\vert^{\frac{1}{2}+\frac{2s}{3}-7\epsilon}}\lesssim 1,\quad \text{ provided that }a<\frac{4s}{3}+2.
\end{equation*}
If $-\frac{1}{2}<s<0$, then
\begin{equation*}
I\lesssim \frac{M^{4\epsilon} \vert k_{3}\vert^{-2s}}{M^{2-a-14\epsilon}\vert k_{1}+k_{3}\vert^{\frac{1}{2}-7\epsilon}}\lesssim \frac{M^{4\epsilon-2s}}{M^{2-a-14\epsilon}}\lesssim 1,\quad \text{ provided that }a<2s+2.
\end{equation*}

\noindent\textbf{Case 4C:} $\vert 2k_{3}^{3}\vert \sim \vert 3(k_{1}+k_{2})(k_{1}+k_{3})(k_{2}+k_{3}) \vert$. Then $\vert k_{3}\vert \sim \vert k_{1}+k_{3}\vert^{\frac{1}{3}} M^{\frac{2}{3}} $. Suppose that $\vert k_{1}+k_{3}\vert\gtrsim  \vert k_{3}\vert$, then $\vert k_{3}\vert^{3}\sim \vert k_{1}+k_{3}\vert M^{2}\gtrsim \vert k_{3}\vert M^{2} $, which yields $\vert k_{3}\vert \gtrsim M $ (contradiction), hence $\vert k_{1}+k_{3}\vert \ll \vert k_{3}\vert $. Therefore, $\vert k_{1}\vert \sim \vert k_{3}\vert $.\par\nopagebreak[4]

\begin{equation*}
I\lesssim \frac{M^{4\epsilon}}{M^{1-a}\vert k_{3}\vert^{2s}}\sim \frac{M^{4\epsilon}}{M^{1+\frac{4s}{3}-a}\vert k_{1}+k_{3}\vert^{\frac{2s}{3}} }.
\end{equation*}
If $s\geq 0$, then
\begin{equation*}
I\lesssim \frac{M^{4\epsilon}}{M^{1+\frac{4s}{3}-a}}\lesssim 1,\quad \text{ provided that }a<\frac{4s}{3}+1.
\end{equation*}
If $s<0$, then

\begin{equation*}
I\lesssim \frac{M^{4\epsilon}\vert k_{1}+k_{3}\vert^{-\frac{2s}{3}}}{M^{1+\frac{4s}{3}-a}}\lesssim \frac{M^{4\epsilon-\frac{2s}{3}}}{M^{1+\frac{4s}{3}-a}}\lesssim 1,\quad \text{ provided that }a<2s+1.
\end{equation*}

\noindent\textbf{Case 5:} $\vert k_{1} \vert,  \vert k_{2} \vert \ll \vert k_{3} \vert \sim M$. Then $\vert k_{4}\vert \sim M $ and $\langle \Phi \rangle \sim M^{3} $.\par\nopagebreak[4]

\begin{equation*}
I\lesssim \frac{\vert k_{1}+k_{2}\vert M^{4\epsilon}}{M^{\frac{7}{2}-a-21\epsilon}\vert k_{1}\vert^{s}\vert k_{2}\vert^{s}}\lesssim \frac{\vert k_{1}\vert^{-s}\vert k_{2}\vert^{1-s}M^{4\epsilon}}{M^{\frac{7}{2}-a-21\epsilon}}.
\end{equation*}
If $s\geq 1$, then 
\begin{equation*}
I\lesssim \frac{M^{4\epsilon}}{M^{\frac{7}{2}-a-21\epsilon}}\lesssim 1,\quad \text{ provided that }a<\frac{7}{2}.
\end{equation*}
If $0\leq s<1$, then
\begin{equation*}
I\lesssim \frac{M^{4\epsilon}\vert k_{2}\vert^{1-s}}{M^{\frac{7}{2}-a-21\epsilon}\vert k_{1}\vert^{s}}\lesssim \frac{M^{4\epsilon+1-s}}{M^{\frac{7}{2}-a-21\epsilon}}\lesssim 1,\quad \text{ provided that }a<s+\frac{5}{2}.
\end{equation*}
If $-\frac{1}{2}<s<0$, then 
\begin{equation*}
I\lesssim \frac{M^{4\epsilon+1-2s}}{M^{\frac{7}{2}-a-21\epsilon}}\lesssim 1, \quad \text{ provided that }a<2s+\frac{5}{2}.
\end{equation*}
\end{proof}
The following propositions concerning the $R_{i}$-terms are established by using arguments analogous to those presented above.

\begin{proposition} \label{prop:R5-R8}
     For $s>-\frac{1}{2} $ and $a< \min\{\frac{3}{2},2s+1,s+1\}$, we have
    \begin{equation*}
\big\Vert R_{5}(u,v,w)+ R_{8}(u,v,w)\big\Vert_{X_{1/6,\delta}^{s+a,b-1}}\lesssim \Vert u\Vert_{X_{-1/6,\delta}^{s,1/2}}\Vert v\Vert_{X_{-1/6,\delta}^{s,1/2}}\Vert w\Vert_{X_{-1/6,\delta}^{s,1/2}}.
\end{equation*}
\end{proposition}
\begin{proposition}\label{prop:R-7}
     For $s>-\frac{1}{2} $ and $a< \min\{1,\frac{4s}{3}+\frac{2}{3}\}$, we have
    \begin{equation*}
\Vert R_{7}(u,v,w)\Vert_{X_{1/6,\delta}^{s+a,b-1}}\lesssim \Vert u\Vert_{X_{1/6,\delta}^{s,1/2}}\Vert v\Vert_{X_{1/6,\delta}^{s,1/2}}\Vert w\Vert_{X_{-1/6,\delta}^{s,1/2}}.
\end{equation*}
\end{proposition}

\begin{proposition}\label{prop:R6-R9-R11}
     For $s>-\frac{1}{2} $ and $a< \min\{1,2s+1\}$, we have
    \begin{equation*}
\big\Vert R_{6}(u,v,w)+ R_{9}(u,v,w)+ R_{11}(u,v,w)\big\Vert_{X_{1/6,\delta}^{s+a,b-1}}\lesssim \Vert u\Vert_{X_{1/6,\delta}^{s,1/2}}\Vert v\Vert_{X_{-1/6,\delta}^{s,1/2}}\Vert w\Vert_{X_{-1/6,\delta}^{s,1/2}}.
\end{equation*}
\end{proposition}

\begin{proposition}\label{prop:R-10}
    For $s>-\frac{1}{2} $ and $a< \min\{2,2s+1\}$, we have
    \begin{equation*}
\Vert R_{10}(u,v,w)\Vert_{X_{1/6,\delta}^{s+a,b-1}}\lesssim \Vert u\Vert_{X_{1/6,\delta}^{s,1/2}}\Vert v\Vert_{X_{1/6,\delta}^{s,1/2}}\Vert w\Vert_{X_{1/6,\delta}^{s,1/2}}.
\end{equation*}
\end{proposition}

\begin{proposition}\label{prop:R-12}
    For $s>-\frac{1}{2} $ and $a<\min \{ 2,2s+1,\frac{4s}{3}+1\}$, we have
    \begin{equation*}
\Vert R_{12}(u,v,w)\Vert_{X_{1/6,\delta}^{s+a,b-1}}\lesssim \Vert u\Vert_{X_{1/6,\delta}^{s,1/2}}\Vert v\Vert_{X_{1/6,\delta}^{s,1/2}}\Vert w\Vert_{X_{-1/6,\delta}^{s,1/2}}.
\end{equation*}
\end{proposition}
\begin{proposition}
   For $s>-\frac{1}{2} $ and $a< \min\{\frac{3}{2},2s+1,s+1\}$, we have 
 \begin{equation*}
\big\Vert \widetilde{R}_{1}(u,v,w)+ \widetilde{R}_{10}(u,v,w)\big\Vert_{X_{-1/6,\delta}^{s+a,b-1}}\lesssim \Vert u\Vert_{X_{1/6,\delta}^{s,1/2}}\Vert v\Vert_{X_{1/6,\delta}^{s,1/2}}\Vert w\Vert_{X_{1/6,\delta}^{s,1/2}}.
\end{equation*}
   
\end{proposition}
\begin{proof}
    The proof is an immediate repetition of the argument used in the proof of Proposition \ref{prop:R5-R8}.
\end{proof}
\begin{proposition}
   For $s>-\frac{1}{2} $ and $a< \min\{1,\frac{4s}{3}+\frac{2}{3}\}$, we have 
    \begin{equation*}
\big\Vert \widetilde{R}_{2}(u,v,w)+ \widetilde{R}_{6}(u,v,w)+ \widetilde{R}_{9}(u,v,w)+ \widetilde{R}_{11}(u,v,w)\big\Vert_{X_{-1/6,\delta}^{s+a,b-1}}\lesssim \Vert u\Vert_{X_{1/6,\delta}^{s,1/2}}\Vert v\Vert_{X_{-1/6,\delta}^{s,1/2}}\Vert w\Vert_{X_{-1/6,\delta}^{s,1/2}}.
\end{equation*}
\end{proposition}
\begin{proof}
  The proof is identical to those of Propositions \ref{prop:R-3}, \ref{prop:R-4}, \ref{prop:R-7} and \ref{prop:R-12} after relabeling the frequency variables according to the permutation $(k_{1},k_{2},k_{3})\rightarrow (k_{2},k_{3},k_{1}) $ and interchanging  frequencies $k_{1} \rightarrow k_{3}$.  
\end{proof}
\begin{proposition}
    For $s>-\frac{1}{2} $ and $a< \min\{1,2s+1\}$, we have
    \begin{equation*}
\big\Vert \widetilde{R}_{3}(u,v,w)+\widetilde{R}_{4}(u,v,w)+\widetilde{R}_{7}(u,v,w)+\widetilde{R}_{12}(u,v,w)\big\Vert_{X_{-1/6,\delta}^{s+a,b-1}}\lesssim \Vert u\Vert_{X_{1/6,\delta}^{s,1/2}}\Vert v\Vert_{X_{1/6,\delta}^{s,1/2}}\Vert w\Vert_{X_{-1/6,\delta}^{s,1/2}}.
\end{equation*}
\end{proposition}

\begin{proof}
  The proof follows by adapting the arguments used in Propositions \ref{prop:R-2} and \ref{prop:R6-R9-R11}   after relabeling the frequency variables according to the permutation
 $(k_{1},k_{2},k_{3})\longrightarrow (k_{3},k_{1},k_{2}) $ and  interchanging frequencies $k_{1} \leftrightarrow k_{3}$.
\end{proof}
\begin{proposition}\label{prop:R-tilde-5-8}
    For $s>-\frac{1}{2} $ and $a< \min\{1,2s+1\}$, we have
    \begin{equation*}
\big \Vert \widetilde{R}_{5}(u,v,w)+\widetilde{R}_{8}(u,v,w)\big\Vert_{X_{-1/6,\delta}^{s+a,b-1}}\lesssim \Vert u\Vert_{X_{-1/6,\delta}^{s,1/2}}\Vert v\Vert_{X_{-1/6,\delta}^{s,1/2}}\Vert w\Vert_{X_{-1/6,\delta}^{s,1/2}}.
\end{equation*}
\end{proposition}
\begin{proof}
   The proof follows exactly the same arguments as the proof of Proposition \ref{prop:R-10} together with the Proposition 3.4 in  \cite{Engin-Burak}.
\end{proof}

\section{Appendix}
We show that the resonance function associated with Proposition \ref{prop:R-3}, in the frequency region $\vert 2k_{3}^{3}\vert \sim \vert 3(k_{1}+k_{2})(k_{1}+k_{3}) (k_{2}+k_{3})\vert$ and $\vert k_{1}\vert ,\vert k_{3}\vert\ll \vert k_{2}\vert \sim M$,
 yields no derivative gain in the denominator. The resonance function corresponding to the $R_{3}$-term  is given by
\begin{equation*}
\Phi =\frac{-2k_{3}^{3}-3(k_{1}+k_{2})(k_{1}+k_{3}) (k_{2}+k_{3})}{6}.
\end{equation*}
Taking $k_{2}=M$ and $k_{1}+k_{3}=1$, we obtain 
\begin{equation*}
\Phi= \frac{-2k_{3}^{3}+3k_{3}^{2}-3k_{3}-3M^{2}-3M}{6}.
\end{equation*}
Under the assumed comparability condition, it follows that $\vert k_{3}\vert \sim M^{2/3}$.
We therefore seek an asymptotic expansion of $k_{3}$ of the form
\begin{equation*}
k_{3}=aM^{2/3}+bM^{1/3}+c+dM^{-1/3}+eM^{-2/3}+fM^{-1}+gM^{-4/3}+...  ,
\end{equation*}
where $a,b,c,d,e,f,g,\ldots$ are real coefficients.
Choosing $a=-\big(\frac{3}{2}\big)^{1/3}$ and $b=0$, and substituting this expansion into the expression for $\Phi$, we obtain
\begin{align*}
\Phi&=(-2a^{3}-3)M^{2}+(3a^{2}-6a^{2}c)M^{4/3}
+(-6a^{2}d-3)M+(-6ac^{2}-6a^{2}e+6ac-3a)M^{2/3}\\
&+(-12acd-6a^{2}f+6ad)M^{1/3}
+(-2c^{3}-12ace-6ad^{2}-6a^{2}g+6ae+3c^{2}-3c)\\
&+\mathcal{O}(M^{-1/3}).
\end{align*}
By successively choosing $c=\frac{1}{2}$, $d=-\frac{1}{2a^{2}}$, $e=-\frac{1}{4a}$ and $f=0$, the coefficients of the positive powers of $M$ vanish. Consequently, the resonance function remains bounded as $M\to \infty$, so that $\langle \Phi \rangle \sim 1$. Thus, the denominator provides no derivative gain in this frequency regime.

\section*{Acknowledgments}
The authors wish to express their sincere gratitude to Mehmet Burak Erdo\u{g}an and O\u{g}uz Y\i lmaz for their valuable discussions and insightful comments.

\bibliographystyle{amsplain}
\bibliography{references}

@article{Xin-Yang-KDV-KDV,
title = {Well-posedness and critical index set of the {C}auchy problem for the coupled {KdV-KdV} systems on $ \mathbb{T} $},
journal = {{D}iscrete and {C}ontinuous {D}ynamical {S}ystems},
volume = {42},
number = {11},
pages = {5167-5199},
year = {2022},
author = {X. Yang and B.-Y. Zhang},
}

@article {Engin-Burak,
    AUTHOR = {Ba\c{s}ako\u{g}lu, E. and G\"urel, T. B.},
     TITLE = {Smoothing and global attractors for the {H}irota-{S}atsuma
              system on the torus},
   JOURNAL = {J. {M}ath. {A}nal. {A}ppl.},
  FJOURNAL = {Journal of {M}athematical {A}nalysis and {A}pplications},
    VOLUME = {525},
      YEAR = {2023},
    NUMBER = {2},
     PAGES = {127244},
}

@ARTICLE{Erdogan-Tzirakis-global-sol-per,
  author={Erdo\u{g}an, M. B. and Tzirakis, N.},
  journal={International {M}athematics {R}esearch {N}otices}, 
  title={Global smoothing for the periodic {K}d{V} evolution}, 
  year={2013},
  volume={20},
  pages={4589-4614},
  }

@article{Bona-Saut-Chen-Media-I,
    author = {J. L. Bona and M. Chen and J. -C. Saut},
    title = {Boussinesq equations and other systems for small-amplitude long waves in nonlinear dispersive media. {I}: {D}erivation and linear theory},
    pages = {283-318},
    journal = {J. {N}onlinear {S}ci.},
    year = {2002},
    Number = {12},  
   
}

@article{Bourgain-partII,
author = {Bourgain, J.},
journal = {Geometric and {F}unctional {A}nalysis},
number = {3},
pages = {209-262},
title = {Fourier transform {R}estriction phenomena for certain lattice subsets and applications to nonlinear evolution equations, Part {II}: {T}he {K}d{V} equation.},
volume = {3},
year = {1993},
}

@article{Bona-Saut-Chen-M-II,
author = {J. L. Bona and M. Chen and J. -C. Saut} ,
year = {2004},
pages = {925-952},
title = {Boussinesq equations and other systems for small-amplitude long waves in nonlinear dispersive media: {II}. {T}he nonlinear theory},
volume = {17},
journal = {J. Nonlinear Sci.}
}

@article {E-compaan,
    AUTHOR = {Compaan, E.},
     TITLE = {Smoothing and global attractors for the {M}ajda-{B}iello
              system on the torus},
   JOURNAL = {Differential and {I}ntegral {E}quations},
  FJOURNAL = {Differential and {I}ntegral {E}quations. {A}n {I}nternational {J}ournal
              for {T}heory \& {A}pplications},
    VOLUME = {29},
      YEAR = {2016},
    NUMBER = {3-4},
     PAGES = {269--308},
      ISSN = {0893-4983}     
}

@article {Oh-diaphantine-cond-in-WPS,
    AUTHOR = {Oh, T.},
     TITLE = {Diophantine conditions in well-posedness theory of coupled
              {K}d{V}-type systems: local theory},
   JOURNAL = {Int. {M}ath. {R}es. {N}ot. },
  FJOURNAL = {International {M}athematics {R}esearch {N}otices. IMRN},
      YEAR = {2009},
    NUMBER = {18},
     PAGES = {3516--3556},   
}

@article{Boussinesq1871,
    author  = {Boussinesq, J. V.},
    title   = {Th{\'e}orie g{\'e}n{\'e}rale des mouvements qui sont propag{\'e}s dans un canal rectangulaire horizontal},
    journal = {Comptes {R}endus de l'{A}cad{\'e}mie des {S}ciences de {P}aris},
    volume  = {73},
    year    = {1871},
    pages   = {256--260}
}

@article{Boussinesq1872,
    author  = {Boussinesq, J. V.},
    title   = {Th{\'e}orie des ondes et des remous qui se propagent le long d'un canal rectangulaire horizontal, en communiquant au liquide contenu dans ce canal des vitesses sensiblement pareilles de la surface au fond},
    journal = {Journal de {M}ath{\'e}matiques {P}ures et {A}ppliqu{\'e}es},
    volume  = {17},
    year    = {1872},
    pages   = {55--108}
}

@article{BabinIlinTiti2011,
    author  = {Babin, A. and Ilin, A. and Titi, E.},
    title   = {On the regularization mechanism for the periodic {Korteweg--de Vries} equation},
    journal = {Communications on {P}ure and {A}pplied {M}athematics},
    volume  = {64},
    number  = {5},
    year    = {2011},
    pages   = {591--648}
}

@article{KenigPonceVega1996,
  author  = {Kenig, C. E. and Ponce, G. and Vega, L.},
  title   = {A bilinear estimate with applications to the {KdV} equation},
  journal = {Journal of the {A}merican {M}athematical {S}ociety},
  volume  = {9},
  number  = {2},
  year    = {1996},
  pages   = {573--603}
}

@article{ChristCollianderTao2003,
  author  = {Christ, M. and Colliander, J. and Tao, T.},
  title   = {Asymptotics, frequency modulation, low regularity
             ill-posedness for canonical defocusing equations},
  journal = {American {J}ournal of {M}athematics},
  volume  = {125},
  year    = {2003},
  pages   = {1235--1293}
}

@article{CollianderEtAl2003,
  author  = {Colliander, J. and Keel, M. and Staffilani, G.
             and Takaoka, H. and Tao, T.},
  title   = {Sharp global well-posedness for {KdV} and modified {KdV}
             on $\mathbb{R}$ and $\mathbb{T}$},
  journal = {Journal of the {A}merican {M}athematical {S}ociety},
  volume  = {16},
  year    = {2003},
  pages   = {705--749}
}

@article{Guo2009,
  author  = {Guo, Z.},
  title   = {Global well-posedness of {Korteweg-de Vries} equation
             in ${H}^{-3/4}(\mathbb{R})$},
  journal = {Journal de {M}ath{\'e}matiques {P}ures et {A}ppliqu{\'e}es},
  volume  = {91},
  number  = {6},
  year    = {2009},
  pages   = {583--597}
}

@article{Kishimoto2009,
  author  = {Kishimoto, N.},
  title   = {Well-posedness of the {C}auchy problem for the
             {Korteweg-de {V}ries} equation at the critical regularity},
  journal = {Differential and {I}ntegral {E}quations},
  volume  = {22},
  number  = {5--6},
  year    = {2009},
  pages   = {447--464}
}

@article{KappelerTopalov2006,
  author  = {Kappeler, T. and Topalov, P.},
  title   = {Global wellposedness of {KdV} in ${H}^{-1}(\mathbb{T}, \mathbb{R})$},
  journal = {Duke {M}athematical {J}ournal},
  volume  = {135},
  number  = {2},
  year    = {2006},
  pages   = {327--360}
}

@article{Molinet2011,
  author  = {Molinet, L.},
  title   = {A note on ill posedness for the {KdV} equation},
  journal = {Differential and {I}ntegral {E}quations},
  volume  = {24},
  number  = {7--8},
  year    = {2011},
  pages   = {759--765}
}

@article{Molinet2012,
  author  = {Molinet, L.},
  title   = {Sharp ill-posedness results for the {KdV} and {mKdV}
             equations on torus},
  journal = {Advances in {M}athematics},
  volume  = {230},
  number  = {4--6},
  year    = {2012},
  pages   = {1895--1930}
}

@article{KillipVisan2019,
  author  = {Killip, R. and Vi\c{s}an, M.},
  title   = {{KdV} is well-posed in ${H}^{-1}$},
  journal = {Annals of {M}athematics},
  volume  = {190},
  number  = {1},
  year    = {2019},
  pages   = {249--305}
}

@article{Ginibre-Velo,
    author = {Ginibre, J. and Tsutsumi, Y. and Velo, G.},
    title = {On the {C}auchy problem for the {Z}akharov system},
    journal = {J. {F}unct. {A}nal.} ,
    volume = {151},
    number = {2},
    year = {1997},
    pages ={384--436}
}

\end{document}